\documentclass[review]{elsarticle}
\usepackage[margin=2cm]{geometry}
\usepackage{hyperref}
\usepackage[]{algorithm2e}
\usepackage{enumitem}

\journal{a journal}

\usepackage{eurosym}
\usepackage{amsmath}
\usepackage{amsfonts}
\usepackage{array}
\usepackage{amssymb}
\usepackage{textcomp}
\usepackage{changepage}
\usepackage{graphicx}
\usepackage{caption}
\usepackage{subcaption}
\usepackage{natbib}
\usepackage{xcolor}
\usepackage{booktabs}
\usepackage{multirow}
\usepackage{epsfig}
\usepackage{setspace}

\biboptions{authoryear}

\begin{document}

\begin{frontmatter}

\title{Analysing the interactions between demand side and supply side investment decisions in an oligopolistic electricity market using a stochastic mixed complementarity problem}


\author{Mel T. Devine\textsuperscript{a,*}, Valentin Bertsch\textsuperscript{b}}
\cortext[mycorrespondingauthor]{Corresponding author}
\ead{mel.devine@ucd.ie}

\address{\textsuperscript{a}Energy Institute, College of Business, University College Dublin, Ireland\\ 
\textsuperscript{b}Chair of Energy Systems and Energy Economics, Ruhr-Universit\"at Bochum, Germany.}

\begin{abstract}
To meet carbon emission targets, governments around the world seek electricity consumers to invest in self-sufficiency technologies such as solar photovoltaic and battery storage. Such behaviour is sought in markets typically characterised by an oligopoly amongst generating firms. In this work, we study the interactions between investment decisions on the demand side and the supply side, and we investigate how price-making behaviour on the supply side affects these interactions compared to a situation with perfect competition. To do so, we introduce a novel stochastic mixed complementarity problem to model several players in an oligopolistic electricity market. On the supply side, we consider generating firms who make operational and investment decisions. On the demand side, we consider both industrial and residential consumers, each of whom may invest in self-sufficiency technologies. The uncertainties of wind and solar generation are the sources of the model’s stochasticity. We apply the model to a case study of a stylised Irish electricity system in 2030. Our results demonstrate that price-making on the supply side increases investment in self-sufficiency on the demand side, leading to a reduction in prices and carbon emissions. We also find that both market power and self-sufficiency alter the investment and decommissioning decisions made by generation firms. Counter-intuitively, we also observe that the absence of a feed-in premium increases investment in solar generation on the demand side. Our findings highlight the importance of including both demand and supply side investment in models of electricity markets characterised by an oligopoly.
\end{abstract}


\end{frontmatter}

\section{Introduction} \label{sec:intro}
To reduce carbon emissions, governments worldwide have developed transformative climate action policies. Such policies typically involve increasing the amount of electricity generated from renewable sources (RES-E), primarily from wind and solar Photovoltaics (PV). Given the ambitious nature of the political targets, significant investment in such sources will be required on both the supply-side (generation) and the demand-side (consumer self-sufficiency). For instance, in its ``Winter Package'' in 2016, the European Commission has stated that it seeks climate-neutrality and wants to put consumers at the centre of the energy transition \citep{european2016clean}. 

Besides this clear political support of so-called prosumers (i.e. consumers that are also producers) by the EU, several energy market developments have contributed to a strong increase in consumer, as well as research interest in, prosumption / consumer self-sufficiency \citep{luthander2015photovoltaic,ellsworth2016conceptualising,bertsch2017drives,schwarz2018self}. On the one hand, learning rates in the growing PV market have led to strong reductions of PV system costs over the past two decades \citep{candelise2013dynamics}. On the other hand, electricity retail prices increased during the same period of time. Consequently, solar PV has reached the so-called ‘grid parity’ in many countries during this time period \citep{khalilpour2015leaving,hagerman2016rooftop,munoz2014power}, i.e. self-consumption of electricity generated from rooftop PV became cheaper than using electricity from the grid at the retail price (and typically also more attractive than feeding electricity from rooftop PV into the grid to earn a feed-in premium). As a result, an increasing number of consumers has started to invest in self-sufficiency in the recent past. The price increase observed during the European energy crisis in 2022 has further increased the profitability of self-sufficiency. 

Such dynamic developments of investment and consumption decisions on the demand side also come along with implications for the supply side. Most importantly, these decisions on the demand side affect the load to be covered by the supply side and hence market prices as well as optimal generation portfolios (including the possibility of decommissioning decisions, in particular of fossil-fuel based generation). Therefore, an in-depth analysis of the interactions between the investment and generation/consumption decisions on the supply and demand side is very important. Hence, in this work, we develop and demonstrate a novel approach to study the optimal generation and storage investment decisions from both the demand and the supply side in addition to decommissioning decisions. Moreover, we seek to study the impacts such decisions have on different types of consumers, generating firms' profits, and carbon emissions.

\subsection{State of the art and research gap}
In order to support decision making in the context of an increasingly complex energy planning, numerous optimisation models have been developed and published in the Operations Research and wider literature \citep{petropoulos2023operational}. These models often adopt linear programming or mixed-integer linear programming approaches and are used to find cost-minimal combinations of energy supply technologies to cover a given energy demand curve under a number of constraints. In the literature, such models are often referred to as energy system models (ESMs). Concerning the analysis of the profitability and implications of consumer self-sufficiency, existing applications of ESMs usually fall into either one of the following three categories: (i) system-scale models that assess the implications of self-sufficiency for the wider energy system, e.g. at a regional, national or continental level \citep[e.g. ][]{schmidt2012regional,merkel2014modelling,fehrenbach2014economic,senkpiel2020systemic,neumann2021costs}, (ii) residential-scale models that mostly focus on household self-sufficiency decisions concerning the investment in or operation of technologies, such as solar PV, battery storage or heat pumps \citep[e.g. ][]{mckenna2017energy,schwarz2018two,schwarz2018self,hayn2018impact,aniello2023shaping,huckebrink2023user,aniello2024subsidies}, or (iii) approaches that soft-link different variants of models in categories (i) and (ii) \citep{gissey2021evaluating,zakeri2021policy,sarfarazi2024}. 

The main limitations of the majority of existing models in the first category include the following: 
\begin{itemize}
    \item First, typical ESMs assume perfect competition and adopt a least-cost central-planning approach. This is critical mainly for two reasons. On the one hand, in today's liberalised energy markets, private-sector investments in (e.g. renewable) energy technologies are required to achieve political targets. To support investors' decisions in renewables, information on the profitability of investment projects is therefore needed. This shortcoming is for instance mitigated by calculating metrics such as the internal rate of return (IRR) or the net present value (NPV) based on ESM results as part of an ex-post analysis \citep[e.g. ][]{finke2023exploring}. On the other hand, in particular in situations where markets become tighter (e.g. as observed during the European energy crisis starting in 2022 or because of policy-driven phase-out decisions and a lack or delay of new investments), concerns around market power may  become more pertinent. For this purpose, game-theoretic equilibrium models have proven useful. There have been several works in the literature that have modelled market power in energy markets \citep{wogrin2013generation,baltensperger2016multiplicity,devine2018shedding, egging2020solving,devine2023strategic,devine2023role}. For a deeper discussion, we refer the reader to \cite{pozo2017basic}. A drawback of most existing game-theoretic equilibrium models, however, is that they usually do not include intricate (technical) details on either the supply or demand sides. Exceptions in the literature include, for example, \cite{ devine2023role}.   
    \item Second, existing models typically model one (usually hourly) market, i.e. one revenue stream for investors. In reality, however, the profitability of investments is often influenced by multiple revenue streams. For instance, to encourage investment in RES-E, several countries have introduced a Feed-in Premium (FiP), i.e. an extra payment - separate to the market price - made to producers of renewable electricity \citep[e.g. ][]{bertsch2017drives,farrell2017specifying}. In addition, an increasing number of countries has introduced a capacity market to encourage investments in firm generation capacity. It is therefore important to study the impact of multiple revenue streams on investment decisions in future energy systems. 
    \item Third, the demand side is often modelled as one system-level aggregated demand curve, i.e. no differences between consumers or consumer groups are considered. On the one hand, this is critical because different consumers and consumer groups have different demand profiles, which affect self-sufficiency levels and hence the economics of self-sufficiency decisions. On the other hand, the use of one aggregated demand curve ignores the fact that elasticities of or costs of reducing energy demand are not uniformly distributed across consumer groups. While game theoretic optimisation approaches usually consider elasticity of energy demand \citep[e.g. ][]{lynch2017investment, newbery2017strategic,  lynch2019impacts, devine2019effect}, the ignorance of differences between consumer groups is a strong simplification \citep{devine2018shedding}. Finally using just one aggregated demand curve also means that interactions between different different types of consumers cannot be analysed. 
    \item Fourth, if different consumers and their decisions are considered explicitly at all, these models usually assume that all consumers participate in the wholesale market directly (which is the only market modelled) and hence pay wholesale prices for their energy demand (\citep[e.g. ][]{lynch2019role}. Again, this is critical since retail prices in reality include several price components in addition to wholesale prices, which has a strong impact on the economics of self-sufficiency decisions. 
    \item Fifth, while ESMs that consider uncertainties (e.g. using stochastic pogramming) do exist \citep[e.g. ][]{ambec2012electricity, moser2016non, wozabal2016effect}, the majority of applications is based on deterministic models \citep{pfenninger2014energy}. This is critical, in particular, for a high-renewable system such as the one examined in this paper, where especially the electricity generation of both wind and PV is highly variable and uncertain.
\end{itemize}

While existing models in the second category consider consumer-specific load profiles and retail prices, including a number of price components on top of wholesale prices, these models typically assume these prices as exogenous inputs \citep[see e.g. ][]{schwarz2018self,aniello2023shaping,nouicer2023economics,aniello2024subsidies}. This means that such models are limited in that no interdependencies between individual consumer decisions and generators, other consumers or the market and market prices more generally can be analysed. In other words, none of these models considers generation investments on both the demand-side and the supply-side and the interactions between them. 

Soft-linking approaches, such as those in the third model category, are generally used increasingly often. Major limitations, however, include that there is usually no general guarantee that the solution obtained will be optimal. Moreover, it is challenging to fully capture the dynamics of the interactions between decisions on the supply and the demand side.

In summary, while some models have been presented in the OR and wider literature that address several of the limitations described above, to the best of our knowledge, there are no existing stochastic game-theoretic equilibrium models with a sufficient level of detail on the supply and demand side that have addressed all of them. While \cite{devine2023role} mitigate quite a few of the limitations (their model considers market power, multiple revenue streams, energy demand profiles and elasticities for different consumer groups, and stochasticity of renewable generation), they do not consider retail prices for different consumer groups or investment decisions on either the supply or demand side.

\subsection{Research questions and contribution}
In the light of the research gap outlined in the previous section, we address the following research questions in this work:

\begin{enumerate}
    \item How does the presence of market power...
    \begin{enumerate}
        \item ... affect (self-sufficiency) investments, consumer tariffs and load shedding  on the demand side?
        \item ... and self-sufficiency investments, in turn affect generation firms’ investment/exit  decisions and profits?
    \end{enumerate}
        \item How does the availability/absence of a FiP as a separate revenue stream affect the above questions?
    \item How do the interactions between market power, investment decisions, and a FiP affect carbon emissions?
\end{enumerate}

To answer these questions, we present a novel stochastic mixed complementarity problem (MCP). MCPs have been used to model various types of energy markets \citep{hobbs2001linear, gabriel2009solving, huppmann2013endogenous, egging2013benders, lynch2015investment, devine2016rolling, devine2023role}. An MCP solves the optimisation problems of several individual players simultaneously and in equilibrium. 

The players we consider on the supply side include power generating firms with different generation portfolios. All firms maximise their profits by optimising the dispatch of their existing assets and making investments into new (renewable and/or conventional) generation under consideration of the RES uncertainty. Firms earn revenues from an energy market, a quantity-based capacity market, and from a FiP. 

On the demand side, we consider a number of different consumer groups as players in our MCP model. We consider commercial/industrial consumers as well as residential consumers. We also distinguish between traditional consumers and prosumers, i.e. consumers that have their own solar PV modules. Prosumers may make further investment decisions in PV and/or storage capacity. Furthermore, we consider an independent battery storage operator who buys and sells electricity from the grid in addition to considering further investments in battery storage.

This MCP builds on the models presented in \cite{lynch2019role} and \cite{devine2023role}. In contrast to \cite{lynch2019role}, however, we consider the impact of market power. We also go beyond \cite{devine2023role} by including generation and storage investment decisions on both the supply side and the demand side (i.e. including self-sufficiency investments) and by considering retail price components for different consumer groups in addition to wholesale prices. 

Thus, to the best of our knowledge, there is no model in the energy market literature such as the one presented in this paper. The novel model features allow us to consider the research questions listed above. 

The remainder of this paper is structured as follows: In section \ref{sec:methodology}, we introduce the mathematical model. In section \ref{sec:data}, we describe the data of our case study. In section \ref{sec:results}, we present our findings. In section \ref{sec:discussion}, we discuss and interpret our findings and summarise the main conclusions.

\section{Methodology} \label{sec:methodology}

\begin{table}[htp!]
\footnotesize
\caption{Indices and sets. }\label{tab:sets}
\begin{tabular}{rl}
\hline
$f\in F$ &Generating firms \\
$t \in T$& Generating technologies\\
$p \in P$& Time periods\\
$k\in K$ & Consumers groups \\
$s \in S$ & Scenarios\\
$h \in H$& Time steps in storage/load shifting period\\
$p' \in P'= \{1,|H|+1,2|H|+1,...\} \subseteq P$& Index representing starting points for storage period\\
 \hline
\multicolumn{2}{l}{Note: sets contain a finite amount of non-zero natural numbers while $|H|$ represents the cardinality of $H$.} 
\end{tabular}%
\end{table}

\begin{table}[htp!]
\footnotesize
\caption{Variables.}\label{tab:variables}
\begin{tabular}{rl}
\multicolumn{2}{l}{Firms' primal variables}\\
\hline
$gen_{f,t,p,s}$& Generation from firm $f$ with technology $t$ in period $p$ and scenario $s$\\
$cap^{\text{bid}}_{f,t}$& Capacity bid of firm $f$ with technology $t$ \\
$inv_{f,t}$&Investment in new generation capacity for firm $f$ with technology $t$\\
$exit_{f,t}$&Decommissioning of old generation capacity for firm $f$ with technology $t$\\
\hline
\multicolumn{2}{l}{Consumers' primal variables}\\
\hline
$g^{\text{ls}}_{k,p,s}$&Load shedding from consumer group $k$ in period $p$ and scenario $s$\\
$g^{\text{up}}_{k,p,s}$& Electricity stored for later time point from consumer group $k$ in period $p$ and scenario $s$\\
$g^{\text{down}}_{k,p,s}$& Electricity used from storage from consumer group $k$ in period $p$ and scenario $s$\\
$g^{\text{pv}}_{k,p,s}$&PV generation from consumer group $k$ in period $p$ and scenario $s$\\
$cap^{\text{STOR}}_{k}$& Consumer group's $k$'s investment in electrical storage \\
$cap^{\text{PV}}_{k}$& Consumer group's $k$'s investment in photovoltaic generation \\
\hline
\multicolumn{2}{l}{Dual variables}\\
\hline
$\gamma_{p,s}$ & System price for time period $p$ and scenario $s$\\
$\kappa$ & Unit capacity price \\
$\lambda^{\#}_{.}$&Lagrange multipliers associated with constraint $\#$ of the firms' problem\\
$\mu^{\#}_{.}$&Lagrange multipliers associated with constraint $\#$ of the firms' problem consumers' problem\\
 \hline
 \multicolumn{2}{l}{Note: '.' is used as a place-holder as the subscripts for both Lagrange multipliers vary depending the on constraint.}
\end{tabular}%
\end{table}

\begin{table}[htp!]
\footnotesize
\caption{Parameters.}\label{tab:parameters}
\begin{tabular}{cp{12.5cm}}
\hline
$PR_{s}$&Probability associated with scenario $s$\\
$MTC_{t}$&Maintenance cost form generating technology $t$\\
$CAP_{f,t}$& Initial generating capacity for firm $f$ with technology $t$\\
$D^{\text{REF}}_{k,p}$& Reference demand for consumer group $k$ in period $p$\\
$LOSS_{k}$& Storage loss factor for consumer group $k$\\
$G_{k}^{\text{LS,MAX}}$& Maximum load shedding for consumer group $k$ in any time period or scenario\\
$INT^{\text{STOR}}_{k}$& Electrical storage/ load shifting capacity for consumer group $k$\\
$INT^{\text{PV}}_{k}$&PV generating capacity for consumer group $k$\\
$FAC_{k}^{\text{STOR}}$& Percentage of electrical storage capacity electricity consumer group $k$ can use in each period and scenario\\
$X^{\text{PREM}}_{k}$& Retail price premium for consumer group $k$\\
$NORM^{\text{G}}_{f,t,p,s}$& Generating profile for firm $f$ with technology $t$ in period $p$ and scenario $s$\\
$NORM^{\text{PV}}_{p,s}$& PV generating profile for period $p$ and scenario $s$\\
$NORM^{\text{WIND}}_{t,p,s}$& Wind profile for wind region (technology) $t$ in period $p$ and scenario $s$\\
$TARGET$ & Capacity target for overall market\\
$X_{t}$ & Feed-In premium for technology $t$\\ 
$X^{\text{PREM}}_{k}$&Retail price premium for consumer group $k$\\
$DR_{t}$ & De-rating factor for technology $t$\\
$A^{.}_{.}$ & Intercept associated with marginal cost functions\\
$B^{.}_{.}$ & Slope associated with marginal cost functions\\
$C^{\text{GEN}}_{t}$ & Marginal cost associated with generating technology $t$\\
$C^{\text{PV}}_{k,p}$& Marginal cost of using PV generation for consumer group $k$ in period $p$\\
$IC^{\text{GEN}}_{t}$&  Annualised investment in generating technology $t$ cost  \\
$IC^{\text{STOR}}_{k}$& Annualised investment in electrical storage cost for consumer group $k$ \\
$IC^{\text{PV}}_{k}$&Annualised investment in PV generating capacity cost for consumer group $k$ \\
 \hline
\end{tabular}%
\end{table}

\begin{table}[htp!]
\footnotesize
\caption{Functions.}\label{tab:functions}
\begin{tabular}{rp{12.5cm}}
\hline
$C^{\text{LS}}_{k,p}(.)$& Load shedding operational cost for consumer group $k$ in period $p$ \\
 \hline
\end{tabular}%
\end{table}

In this section, the methodology is described. We develop a stochastic MCP to represent an electricity market with two types of players: generation firms and electricity consumer groups. 

Firms receive revenues from energy and capacity markets as well as a FiP and seek to maximise their profits. They may hold multiple generating units with the technologies considered being coal, oil, Combined Cycle Gas Turbine (CCGT), hydroelectric, wind and solar. In Section \ref{sec:data} we consider different wind regions of Ireland with varying, but correlated, hourly capacity factors. Wind for the different regions are represented in the model as different technologies, i.e., wind from one region is considered as a different technology from wind in another region. Firms are distinguished by the initial generating portfolio they hold but may also invest in additional capacity in any of these technologies. All firms may be modelled as either price-takers or price-makers.

The capacity payment mechanism is a quantity based mechanism. Under such a framework, a policy maker, regulator or Transmission System Operator (TSO), exogenously selects a quantity of capacity required for a given period (eg., one year). An auction is held wherein generators submit capacity bids and the total amount sold equal the MW target. Generation firms compete in the auction to hold the options and therefore gain a fixed sum of money per the time period considered to compensate them for each installed unit of capacity.

Consumers seek to minimise the cost of meeting their demand. They may generally do so by utilising a range of possible demand-side flexibility measures, such as load shedding, load shifting or electrical storage, or PV generation. We do not model individual consumers but rather consider different consumer groups whose decisions represent the aggregate actions of consumers in these groups. Consumer groups are distinguished by whether they represent industrial/commercial or residential consumers and whether or not they can invest in additional demand-side flexibility capacity (storage and PV). 

The stochasticity of the model arises from the uncertainty surrounding wind and PV power. Thus, each scenario in our model represents different RES generation profiles, i.e. varying levels of wind and solar power availability at each point in time. Each of the generation firms and consumer groups considered have separate optimisation problems that are connected through market clearing conditions. The stochastic MCP is made up of these market clearing conditions along with the Karush-Kuhn-Tucker (KKT) conditions for optimality from each of the players. Thus, the MCP solves the optimisation problem of each player simultaneously and in equilibrium. The KKT conditions are presented in the Online Appendix\footnote{\url{https://figshare.com/articles/journal_contribution/Online_Appendix_for_Devine_and_Bertsch_2024/25764498}}.

Throughout this section the following conventions are used: lower-case Roman letters indicate indices or primal variables, upper-case Roman letters represent parameters (i.e., data, functions), while Greek letters indicate prices or dual variables. The variables in parentheses alongside each constraint in this section are the Lagrange multipliers associated with those constraints.

\subsection{Firm $f$'s problem} \label{subsec:firms_problem}
Firm $f$ maximises its expected profits (revenues less cost) by choosing the amount of generation ($gen_{f,t,p,s}$), capacity bid ($cap^{\text{bid}}_{f,t}$), investment in new capacity ($inv_{f,t}$) and decommissioning of existing capacity ($exit_{f,t}$). It considers revenues received from a capacity and an energy market as well as a FiP for RES generation. Its costs include generation costs and investment costs in addition to costs incurred for maintaining all of its units. Firm $f$'s generation decisions are scenario dependent while all of its other decision variables are scenario independent. Firm $f$'s optimisation problem is:

\begin{subequations}
\begin{equation}\label{eqn:gen_obj}
\begin{split}
\max_{\substack{ gen_{f,t,p,s}, cap^{\text{bid}}_{f,t},\\ 
inv_{f,t}, exit_{f,t}\\
 }}\>\>
&\sum_{t,p,s}\bigg( PR_{s}\times gen_{f, t,p,s} \times \big (  \gamma_{p,s}+X_{t}   -  C^{\text{GEN}}_{t} \big )\bigg )\\
-&\sum_{t} \bigg(   IC^{\text{GEN}}_{t}\times inv_{f,t}   +\big(inv_{f,t}+CAP_{f,t} - exit_{f,t}\big)\times MTC_{t}  \bigg)\\
+&\sum_{t} DR_{t}\times \kappa \times cap^{\text{bid}}_{f,t},\\
\end{split}
\end{equation}
subject to:
\begin{equation}\label{eqn:gen_energy_con}
 gen_{f, t,p,s} \leq (CAP_{f,t}+inv_{f,t} - exit_{f,t})\times NORM^{\text{G}}_{f,t,p,s}, \>\> \forall t,p,s, \>\> (\lambda^{1}_{f,t,p,s}),
\end{equation}
\begin{equation}\label{eqn:gen_cap_bid_con}
 cap^{\text{bid}}_{f,t} \leq CAP_{f,t}+inv_{f,t} - exit_{f,t}, \>\> \forall t, \>\> (\lambda^{2}_{f,t}),
\end{equation}
\end{subequations}
where $t$ represents the different generating technologies, $p$ represents timesteps and $s$ represents different wind and solar scenarios of the stochastic optimisation problem. The parameters $IC_{t}$ and $MTC_{t}$ are cost of investment and maintenance for technology $t$ respectively. The marginal cost of generating with technology $t$ is $C^{\text{GEN}}_{t}$.

The parameters $X_{t}$, $PR_{s}$, and $DR_{t}$ represent the Feed-in Premium associated with technology $t$, the probability associated with scenario $s$, and the de-rating factor associated with technology $t$, respectively. In this work a de-rating factor reflects the proportion of its overall capacity a technology can provide to meet the capacity target. 

The parameter $CAP_{f,t}$ is the initial installed capacity firm $f$ has for technology $t$ while $NORM^{\text{G}}_{f,t,p,s}$ represents the availability for technology $t$ in timestep $p$ and scenario $s$ for firm $f$. We assume $NORM^{\text{G}}_{f,t,p,s}=1$ $\forall f,t,p,s$ for conventional generation. For wind and PV generation, $NORM^{\text{G}}_{f,t,p,s}=NORM^{\text{WIND}}_{t,p,s}$ and $NORM^{\text{G}}_{f,t,p,s}=NORM^{\text{PV}}_{p,s}$, respectively. The values for $NORM^{\text{WIND}}_{t,p,s}$ and $NORM^{\text{PV}}_{p,s}$ can be found in the online appendix\footnotemark[1]. They take a value between zero and one depending on the timestep and scenario representing the intermittency and uncertainty of RES generation. 

Constraint \eqref{eqn:gen_energy_con} ensures that firm $f$ with technology $t$ cannot generate more than its capacity times its capacity factor, while constraint \eqref{eqn:gen_cap_bid_con} ensures it cannot make a capacity bid greater than the installed capacity. Each of firm $f$'s primal (decision) variables are also constrained to be non-negative.

The capacity price paid for each unit of capacity accepted is $\kappa$. It is exogenous to firm $f$'s problem but is a variable of the overall problem, determined via the market clearing condition \eqref{eqn:MCC_cap}. The energy price at each period and scenario is $\gamma_{p,s}$. If firm $f$ is a price taker then its decision variables cannot affect this price (i.e., $\frac{\partial \gamma_{p,s}}{\partial gen_{f, t,p,s}}=0$). In this case, $\gamma_{p,s}$ is exogenous to firm $f$'s problem but is a variable of the overall problem, determined via the market clearing condition \eqref{eqn:MCC_demand}. If firm $f$'s generating unit for technology $t$ is assumed to be a price-maker unit, then its generation decision variable for that unit ($gen_{f,t,p,s}$) can affect the energy price. As a result, we derive the following relationship between the energy price and generation:
\begin{equation}\label{eqn:gen_price_rel}
\begin{split}
gen_{f,t,p,s} = &\sum_{k} \big(D^{\text{REF}}_{k,p}  +  g^{\text{up}}_{k,p,s} -    g^{\text{down}}_{k,p,s}- g^{\text{pv}}_{k,p,s}\big)
-\sum_{k}\frac{\gamma_{p,s}-A^{\text{LS}}_{k,p}+\frac{1}{PR_{s}}(\mu^{1}_{k,p,s}+\mu^{8}_{k,p,s}-\mu^{1A}_{k,p,s})}{2B^{\text{LS}}_{k,p}} 
\end{split}
\end{equation}
where the parameters and variables not mentioned already are parameters and variables from the consumers' problem (section \ref{subsec:consumers_problem}) and hence exogenous to firm $f$'s problem. Equation \eqref{eqn:gen_price_rel} is determined by combining market clearing condition \eqref{eqn:MCC_demand} with the KKT conditions\footnotemark[1] that determine how consumers shed their load (A.8) and the KKT conditions that determine how firms generate (A.1) (all except firm $f$'s unit $t$). The remaining KKT conditions cannot be used as they cannot be substituted into \eqref{eqn:MCC_demand}. For price-making firms, this relationship is substituted into firm $f$'s objective function \eqref{eqn:gen_obj} leading to 
\begin{equation}\label{eqn:gen_price_gen_rel}
\frac{\partial \gamma_{p,s}}{\partial gen_{f, t,p,s}} = -\bigg(\frac{1}{\sum_{k}\frac{1}{2B^{\text{LS}}_{k,p}}}\bigg).
\end{equation}
Equation \eqref{eqn:gen_price_gen_rel} represents how firm $f$ believes it can influence $\gamma_{p,s}$ with its generation decisions. Consequently, equation \eqref{eqn:gen_price_gen_rel} may over or underestimate firm $f$'s level of market power. As outlined in \cite{devine2023role}, making such assumptions around how price-making firms believe they can influence the market price is a limitation of the modelling approach considered in this work. However, modelling firms' beliefs on how they can influence the market price, also known as conjectural variations, is a practise that has been widely used in the energy market modelling literature \citep{egging2008complementarity, haftendorn2010modeling, huppmann2014market, baltensperger2016multiplicity, egging2016risks}.

If firm $f$ is a price-taker, then its problem is convex. If firm $f$'s generating unit for technology $t$ is a price-maker, then its problem is strictly convex, assuming $B^{\text{LS}}_{k,p} > 0$, $\forall k,p,t$. The firms' KKT conditions are presented in the online appendix\footnotemark[1]. Note: $PM_{f,t}$ is a binary parameter that is used in the KKT conditions to indicate whether firm $f$'s generating unit for technology $t$ is a price making unit ($PM_{f,t}=1$) or price taking unit ($PM_{f,t}=0$).

\subsection{Consumer group $k$'s problem} \label{subsec:consumers_problem}
Consumer group $k$ seeks to  minimise the cost of meeting their expected demand by choosing from a  range of demand side flexibility measures:
\begin{enumerate}
\item \textit{Load shedding}:  consumer group $k$ may reduce  their demand at any period $p$ or scenario $s$ The decision variable that represents the amount they do so by is $g^{\text{LS}}_{k,p,s}$.
\item \textit{Electrical storage}: instead of obtaining electricity from the grid, consumers may obtain electricity from electrical storage. However, they must increase their consumption in previous timestep(s) to ensure there is availability of electricity in storage.  Storing electricity and using electricity to meet demand in period $p$ and scenario $s$ are represented by the decision variables $g^{\text{UP}}_{k,p,s}$ and $g^{\text{DOWN}}_{k,p,s}$, respectively. Consumer group $k$ may also increase their electrical storage capacity which is represented by the investment decision variable $cap^{\text{STOR}}_{k}$.
\item \textit{Photovoltaic (PV) generation}: instead of obtaining electricity from the grid, consumers may also obtain electricity from their own PV-generation unit. The decision variable that represents the amount they do so by is $g^{\text{PV}}_{k,p,s}$ while $cap^{\text{PV}}_{k}$ represents their investment in PV generation decision variable. The availability of PV power varies across timesteps and scenarios.
\end{enumerate}

Consumer group $k$'s investment decisions are scenario independent while all their other decision variables are scenario dependent.
Consumer groups $k$'s reference demand for period $p$ is $D^{\text{REF}}_{k,p}$. This is the demand consumer group $k$ would have in the absence of any demand side flexibility measures. Consumer group $k$'s demand in period $p$ and scenario $s$ is its reference demand plus any storage charging less any load shedding, PV generation or storage discharging. The cost of meeting a unit of this demand for consumer group $k$ is the wholesale price $\gamma_{p,s}$ plus the retail price premium $X^{\text{PREM}}_{k}$.

Consumer group $k$'s optimisation problem is:
\begin{subequations}\label{eqn:con_problem}
\begin{adjustwidth}{-1.5cm}{-1.5cm}
\begin{equation}\label{eqn:obj_act}
\begin{split}
\min_{\substack{  g^{\text{LS}}_{k,p,s},g^{\text{UP}}_{k,p,s},g^{\text{DOWN}}_{k,p,s},\\ g^{\text{PV}}_{k,p,s} cap^{\text{STOR}}_{k}, cap^{\text{PV}}_{k} } }&
\sum_{s,p} PR_{s}  \bigg( (\gamma_{p,s}+X^{\text{PREM}}_{k} )\times \big(D^{\text{REF}}_{k,p} -   g^{\text{LS}}_{k,p,s}   +  g^{\text{UP}}_{k,p,s} -   (1-LOSS_{k})g^{\text{DOWN}}_{k,p,s}- g^{\text{PV}}_{k,p,s}\big)\\ &  +
g^{\text{LS}}_{k,p,s}\times C^{\text{LS}}_{k,p}(  g^{\text{LS}}_{k,p,s}) + 
g^{\text{PV}}_{k,p,s}\times C^{\text{PV}}_{k,p} \bigg)
   +IC^{\text{STOR}} \times cap^{\text{STOR}}_{k} +IC^{\text{PV}}_{k} \times cap^{\text{PV}}_{k},
\end{split}
\end{equation}
\end{adjustwidth}
subject to
\begin{eqnarray}
g^{\text{LS}}_{k,p,s} &\leq &   G_{k}^{\text{LS,MAX}},\>\> \forall p,s, \>\> (\mu^{1}_{k,p,s}),\label{eqn:LS_max_con}\\
g^{\text{LS}}_{k,p,s} &\geq &   0,\>\> \forall p,s, \>\> (\mu^{1\text{A}}_{k,p,s}),\label{eqn:LS_non_neg_con}\\
g^{\text{UP}}_{k,p,s} &\leq &   FAC_{k}^{\text{STOR}}\times(INT^{\text{STOR}}_{k}+cap^{\text{STOR}}_{k}),\>\> \forall p,s, \>\> (\mu^{2}_{k,p,s}),\label{eqn:UP_max_con}\\
g^{\text{DOWN}}_{k,p,s} &\leq &   FAC_{k}^{\text{STOR}}\times(INT^{\text{STOR}}_{k}+cap^{\text{STOR}}_{k}),\>\> \forall p,s, \>\> (\mu^{3}_{k,p,s}),\label{eqn:DOWN_max_con}\\
g^{\text{PV}}_{k,p,s} &\leq & NORM^{\text{PV}}_{p,s} \times (INT^{\text{PV}}_{k}+cap^{\text{PV}}_{k}) ,\>\> \forall p,s, \>\> (\mu^{5}_{k,p,s}),\label{eqn:PV_max_con}\\
\sum^{p'+h-1}_{e=p'}\big( g^{\text{UP}}_{k,e,s}-g^{\text{DOWN}}_{k,e,s} \big) &\leq & INT^{\text{STOR}}_{k}+cap^{\text{STOR}}_{k},\>\>\forall s,p',h ,\>\> (\mu^{6}_{k,p',h,s}),\label{eqn:con_stor_upper}\\
\sum^{p'+h-1}_{e=p'}\big( g^{\text{DOWN}}_{k,e,s}-g^{\text{UP}}_{k,e,s} \big) &\leq & 0,\>\>\forall s,p',h ,\>\> (\mu^{7}_{k,p',h,s}),\label{eqn:con_stor_lower}\\
g^{\text{LS}}_{k,p,s}    +   (1-LOSS_{k})g^{\text{DOWN}}_{k,p,s}+g^{\text{PV}}_{k,p,s}   &\leq&  D^{\text{REF}}_{k,p}+  g^{\text{UP}}_{k,p,s},  \>\> \forall p,s, \>\> (\mu^{8}_{k,p,s}), \label{eqn:con_full_market}
\end{eqnarray}
\end{subequations}
where the parameter $LOSS_{k}$ is a percentage that represents the loss factor associated with storing electricity for consumer group $k$. 

The marginal cost function associated with load shedding is
\begin{eqnarray}
C^{\text{LS}}_{k,p}(x)&=& A^{\text{LS}}_{k,p}+ B^{\text{LS}}_{k,p}x,\label{eqn:shedding_cost}
\end{eqnarray}
while $C^{\text{PV}}_{k,p}$ represents the marginal of PV generation.

The parameter $G_{k}^{\text{LS,MAX}}$ is the maximum amount of electricity that consumer group $k$ can shed their load by in each period and scenario while $INT^{\text{PV}}_{k}$ represent the initial PV generation capacity consumer group $k$ holds. Furthermore, the parameter $INT^{\text{STOR}}_{k}$ represents the initial amount electricity consumer group $k$ may keep in storage while $FAC_{k}^{\text{STOR}}$ represents the percentage that consumer group $k$ can use of their storage in each period and scenario $s$. 

Constraints \eqref{eqn:LS_max_con} - \eqref{eqn:PV_max_con} constrain the amount of electricity consumer group $k$ can generate or store in time period $p$ and scenario $s$. Constraint \eqref{eqn:con_stor_upper} ensures consumer group $k$ cannot, over a $|H|$-timestep period, store more electricity than its storage capacity The index $h$ presents timesteps in a $|H|$-timestep period while the index $p'$ represents starting points for the storage period, i.e., $p' \in P'= \{1,|H|+1,2|H|+1,...\} \subseteq P$.  Constraint \eqref{eqn:con_stor_lower} ensures consumer group $k$ cannot, over the same $|H|$-step time period, use more electricity for meeting demand than what has already been stored. Constraints \eqref{eqn:con_stor_upper} and \eqref{eqn:con_stor_lower}  also ensure that any electricity stored in a $|H|$-timestep period cannot be used in any other $|H|$-timestep period. In reality, if a consumer stores electricity, they will be able to use it in any time period in the future. However, because  we set $|H|=48$ in Section \ref{sec:data}, we believe this simplification is reasonable given the daily trough and peak structure of electricity demand. In addition, this simplification reduces the size and complexity of the model and hence increases computational efficiency.

Constraint \eqref{eqn:con_full_market} ensures any electricity consumer group $k$ generate on their own, either from PV, or storage must be less than their reference demand plus any increased consumption due to storing electricity. In other words, consumer group $k$'s own generation cannot be  used to meet other consumers' demand. All of consumer group $k$'s primal (decision) variables are also constrained to be non-negative. 

Consumer group $k$'s problem is convex, assuming all values for $B^{\text{LS}}_{k,p}$ are non-negative. Its KKT conditions are presented in the online appendix\footnotemark[1]. In Section \ref{sec:data} and \ref{sec:results} we also an independent storage operator. Its optimisation problem is a reduced form of a consumer group's problem.

\subsection{Market clearing conditions} \label{subsec:clearing}
The $|F|+|K|$ optimisation problems are connected via the following market clearing conditions:
\begin{subequations}\label{eqn:MCCs}
\begin{equation}\label{eqn:MCC_demand}
\sum_{f,t} gen_{f,t,p,s} = \sum_{k} \big(D^{\text{REF}}_{k,p} -   g^{\text{ls}}_{k,p,s}   +  g^{\text{up}}_{k,p,s} -    (1-LOSS_{k})g^{\text{down}}_{k,p,s}- g^{\text{pv}}_{k,p,s}\big),\>\> \forall p,s,\>\> (\gamma_{p,s}),
\end{equation}
\begin{equation}\label{eqn:MCC_cap}
\sum_{f,t} DR_{t}\times cap^{\text{bid}}_{f,t} = TARGET,\>\>  (\kappa),
\end{equation}
\end{subequations}
Market clearing condition \eqref{eqn:MCC_demand} ensures that the total amount of electricity generated by the firms must equal the sum of the consumers' demand. Consumers' demand consists of their reference demand plus any electricity they store less any electricity they shed or generate themselves. Market clearing condition \eqref{eqn:MCC_cap} ensures that the capacity bids of the firms must equal the capacity target level. The variables $\gamma_{p,s}$ and $\kappa$  are the prices and free Lagrange multipliers associated with conditions \eqref{eqn:MCC_demand} and \eqref{eqn:MCC_cap}, respectively.

As each of the optimisation problems are convex, the KKT conditions are both necessary and sufficient for optimality for each type of player \citep{gabriel2012complementarity}. Thus, the stochastic MCP consists of conditions (A.1) - (A.21) (see the Online Appendix\footnotemark[1]) in addition to the market clearing conditions \eqref{eqn:MCCs}.

\subsection{Solving the problem}\label{subsec:solving}

For computational efficiency, the MCP is solved using Benders Decomposition. The Benders Decomposition algorithm used in this paper is described in detail in Online Appendix B\footnotemark[1]. It follows those used to solve MCPs in \cite{egging2013benders}, \cite{gabriel2010benders} and \cite{bertsch2018analysing}; the overall scheme is presented in Online Appendix B.1. Solving the stochastic MCP in this paper involves two steps:

\begin{enumerate}
\item The MCP is solved for a selection of 24 days in hourly resolution. These days represent 8-day periods in winter, spring/autumn and summer, resulting in a total of 24 $\times$ 8 $\times$ 3 = 576 hourly time steps. We choose 8-day periods, instead of 7-day periods (weeks), to ensure an even number of days in each period. This is because we choose $H=48$ timesteps for each storage/load shifting period and thus require an even number of days to ensure there is equal amount of timesteps in each storage/load shifting period. The revenues and costs from each hour in the objective function of each player are multiplied by weighting factors of 11.375, 22.750 or 11.625 depending on whether that hour represents a week in winter, spring/autumn or summer, respectively. This ensures that the 576 hourly time steps represent a full year of 8784 hours. To determine the weighting factors we divided the whole year into seasons and assigned 91 days (13 weeks) to winter, spring and autumn respectively and 93 days to summer. To reduce the computational intensity of the model, spring and autumn are represented by the same week explaining the higher weighting factor of the time steps representing these seasons. 

For this first step, a Benders Decomposition algorithm is used to solve the MCP. Benders Decomposition is a solution algorithm that has been shown to solve stochastic MCPs in a computational efficent manner \citep{egging2013benders, gabriel2010benders}. Instead of solving for all variables simultaneously it divides the MCP into multiple smaller MCPs and solves the overall problem iteratively. At each iteration, a first-stage master problem MCP is solved for a subset of variables. The values of the remaining variables are determined by second-stage sub-problem MCPs.

In the model presented above, the different timesteps and scenarios are connected through the investment, exit and capacity bid variables. If these primal variables are fixed to specific
values, the model could be solved for each storage period and scenario separately. Hence these variables are the first-stage variables and are also known as the complicating variables. 

The master problem MCP of this work is presented in Online Appendix B.2. Its initial constraints are those of the overall problem that involves these variables. Constraints not involving any of these variables are not included in the master problem MCP. At each iteration, new constraints, known as Benders cuts, are added to the master problem MCP using the results from the preceding sub-problems. Consequently, the solution obtained is improved at each iteration. The algorithm stops when the first-stage variables produced by the master problem stop changing. A metric derived in \cite{gabriel2010benders} is used to measure this convergence; see Online Appendix B.3.

The sub-problem MCPs can be represented by the overall MCP presented in Section \ref{sec:methodology}. However, each sub-problem is solved with investment, exit and capacity bid variables fixed at the values determined by the master problem MCP of the same iteration. The second-stage variables are all scenario specific and, in this work, are the firms' generation variables in addition to each of the consumers' load shedding, PV, charging, and discharging decision variables.

\item The optimal investment, exit and capacity bid variables from the first step are fixed as parameters in the second step. The MCP is then solved 93 times (without Benders decomposition), each time representing a different 48-hour period in a 366-day year. As before, we choose a 366-day year to ensure an even number of days and thus ensure there is an equal amount of timesteps in each 48-hour storage period. As inter-temporal constraints do not exceed beyond 48 hours, splitting the MCP into smaller problems is equivalent to solving a single MCP with 8784 timesteps, as investment, exit, and capacity bid decisions are fixed. However, splitting the problem up into multiple smaller problems is more computationally efficient \citep{devine2016rolling}. In addition, we believe it is reasonable to assume that investment/exit decisions are taken at times separate to generation, storage and load shedding decisions.
\end{enumerate}

The outputs of the model are the optimal investment, exit and capacity bid decisions from the first step and the optimal generation, load shedding and storage decisions and resulting prices from the second step. From this model output, we can calculate consumer costs, generators profits and CO$_2$ emissions. 

\section{Input data}\label{sec:data}
The Irish government in its climate action plan \citep{CAP2023} has committed to significant RES-E deployment on both the demand and supply sides. We apply the model described in section \ref{sec:methodology} to a case study based on the future Irish power system. For this purpose, we mainly use data for 2030 from \citet{eirgrid2019}. The equilibrium problem is solved over $|P|=8784$ hourly time periods. Demand side data are described in section \ref{sec:data_demand} while in section \ref{sec:data_SO} we describe data for an independent storage operator. We describe input data related to the supply side in section \ref{sec:data_conventional} and data related to renewable generation in section \ref{sec:input-data}. 

\subsection{Demand side data}
\label{sec:data_demand}
We consider four different consumer groups on the demand side. These include commercial/industrial ($k=1,2$) as well as residential consumers ($k=3,4$) and within each of these, we distinguish between traditional consumers ($k=1,3$) and prosumers ($k=2,4$), i.e. consumers that have their own generation from solar PV modules and/or invest in storage capacity. Traditional consumer groups' investment decisions are fixed at zero.

In Section \ref{sec:results} we vary the percentage of consumers who are prosumers by varying the size of the different consumer groups' reference demand. For instance, if the residential prosumer group's reference demand is $\alpha\%$ of total residential reference demand for period $p$ and scenario $s$, then the reference demand of residential consumer group that cannot invest in demand-side flexibility is $(1-\alpha)\%$ of total residential reference demand for the same period. Likewise for industrial consumers. 

Total reference demand profiles of the industrial and residential consumer groups can be found in the Online Appendix\footnotemark[1]. Based on \cite{eirgrid2019}, we consider an overall annual electricity demand of 39 TWh and a peak demand of 6503 MW. We calculate the capacity target to be equal to the system peak demand.

We allow consumers to shed all their load. The parameters for the marginal load shedding cost functions (\ref{eqn:shedding_cost}) follow from from \citet{devine2018shedding} and can be found in the Online Appendix\footnotemark[1]. \citet{devine2018shedding} use value of lost load (VOLL) estimates by \citet{leahy2011estimate} for different types of consumers in Ireland to derive load shedding cost functions for industrial and residential consumers. 

All four consumer group are assumed to have no initial storage capacity. For consumer groups who cannot invest in storage we fix $cap^{\text{STOR}_{k=1,3}}=0$ while for prosumers, the storage investment costs is $IC^{\text{STOR}}=$\euro 32,248/MW. This assumes a lifespan of 20 years and an interest rate of 7.1\% \citep{pietzcker2021tightening}. We assume that consumers may charge (or discharge) the full storage capacity in any given hour, i.e., $FAC^{\text{STOR}_{k=2,3}}=1$. The storage loss parameter is 10\% \citep{pietzcker2021tightening}, i.e., $LOSS_{k=2,4}=0.1$.

The prosuming residential and industrial consumer groups each have initial PV capacity of 62.5MW \citep{eirgrid2019}. We assume industrial prosumers invest in large-scale PV generation while residential prosumers invest in small-scale PV generation. That is, we assume annualised PV investment costs are $IC^{\text{PV}}_{k=2}=$ \euro 44,119/MW and $IC^{\text{PV}}_{k=4}=$ \euro 61,913/MW. These assumes a lifespan of 25 years and interest rate 0f 7.1\% and 6.8\%, respectively, \citep{pietzcker2021tightening,fraunhofer2018stromgestehungskosten}. Installation costs for residential PV are higher due to the economics of scale associated with industrial PV.
 
We assume that the retail price premium $X^{\text{PREM}}_{k=1,2}=$ \euro 35/MWh for industrial consumers and $X^{\text{PREM}}_{k=3,4}=$ \euro 125/MWh for residential consumers. These were obtained by comparing average Irish wholesale prices\footnote{\url{https://www.statista.com/statistics/1271371/ireland-monthly-wholesale-electricity-price/}} for 2022 with average industrial and residential prices\footnote{\url{ https://www.seai.ie/data-and-insights/seai-statistics/key-statistics/prices/}} in Ireland for 2022 and 2021.

\subsection{Independent storage operator}\label{sec:data_SO}
In this work, we consider an independent Storage Operator (SO) as a player in the equilibrium model. The SO's optimisation problem is a reduced form of a consumer group's problem (equations \eqref{eqn:con_problem}), labelled $k=5$. The SO buys electricity to charge its batteries and then subsequently discharges the batteries and sells that electricity. It does so to minimise its costs (maximise its profits). It buys and sells electricity at the wholesale price and thus $X^{\text{PREM}}_{k=5}=0$. 

The SO does not have any demand ($D^{\text{REF}}_{k=5,p,s}=0 \> \forall p,s$). Consequently, it cannot invest in, nor generate, PV generation, nor can it shed any load. To allow the SO to sell its discharged electricity to the grid, constraint \eqref{eqn:con_full_market} is excluded from its optimisation problem. 

In line with \cite{eirgrid2019} values for 2030, the SO has an initial capacity of $CAP^{\text{STOR}}_{k=5}=700$MW and may invest in further storage capacity at the same cost of the consumers (Section \ref{sec:data_demand}). Furthermore, as with the consumers, it may charge (discharge) its capacity fully in any given hourly period and we assume a storage loss factor of 10\%.

\subsection{Supply side data}
\label{sec:data_conventional}
On the supply side, we consider $|F|=5$ generation firms who hold and may invest in multiple generating technologies. The technologies considered are coal, CCGT, oil, hydroelectric, wind and solar. For wind, we consider three different investment options (wind1, wind2, and wind3), each representing different wind regions in Ireland and hence generating profiles; see Section \ref{sec:input-data}, \cite{bertsch2018analysing} and \cite{devine2023role} for further details. 

The initial installed capacity portfolio is presented in Table \ref{tab:firm_data}. The maximum capacity values for each technology are broadly based on \cite{eirgrid2019}.

\begin{table}[htp!]
 \centering
\footnotesize
\caption{Initial power generation portfolio by firm ($CAP_{f,t}$). }\label{tab:firm_data}
\begin{tabular}{p{4cm}p{1.5cm}p{1.5cm}p{1.5cm}p{1.5cm}p{1.5cm}}
\hline
Technology & Firm 1 & Firm 2 & Firm 3 & Firm 4 & Firm 5\\
\hline
Coal&570&0&476&0&0\\
Oil&338&696&258&0&0\\
CCGT&1992&464&648&1370&1370\\
Hydro&508&0&0&0&0\\
Wind1&709&241&288&285&285\\
Wind2&882&300&358&355&355\\
Wind3&605&206&245&243&243\\
Solar&49&17&20&20&20\\
\hline
\end{tabular}%
\end{table}

Table \ref{tab:data_gen} shows the costs and carbon emissions associated with the different generating technologies. To calculate marginal costs for coal, oil and CCGT, we assume power plant efficiencies 38\%, 42\% and 60\%, respectively, and a carbon price of \euro 129/t CO$_{2}$. Moreover, we assume the fuels associated with these technologies cost \euro 3/GJ, \euro 14.3/GJ, and \euro 7.1/GJ, respectively. 

To calculated annualised investment costs, we assume investment costs of \euro 1800/kW, \euro 400/kW, \euro 900/kW, \euro 2500/kW, \euro 1137/kW, and \euro 395/kW and lifespans of 45, 45, 40, 80, 25, and 25 years for coal, oil, CCGT, hydro, wind, and solar, respectively. We also assume maintenance costs are 2\%, 4\%, 3\%, 3\%, 3\%, and 2\%  of fixed costs for each of these technologies, respectively. All these parameters are based off values for 2030 from \cite{pietzcker2021tightening}. 

We assume interest rates of 10.2\% for thermal generation, 7.5\% for wind generation and 7.1\% for solar generation \citep{fraunhofer2018stromgestehungskosten}. We assume a FiP ($X_{t}$) of \euro 23/MW \citep{farrell2017specifying} for wind and solar and a FiP of zero for all other technologies (In Section \ref{sec:results} we also consider cases where $X_{t}=0$ $\forall t$. ). 

\begin{table}[htp!]
 \centering
\footnotesize
\caption{Summary of techno-economic input data of considered supply side technologies. }\label{tab:data_gen}
\begin{tabular}{p{3cm}p{2.25cm}p{2.25cm}p{2.5cm}p{2.5cm}}
\hline
Technology & Annuity of specific invest  & Fixed O\& M costs  & Marginal power gen. costs at intercept & Spec. CO$_2$ emissions \\
&($IC^{\text{GEN}}_{t})$ &($MTC_{t}$)&($C^{\text{GEN}}_{t}$)&-\\ 
&(\euro /GW y) &(\euro/GW y)&(\euro/MWh$_{el}$)&(t CO$_2$/MWh$_{el}$)\\
\hline
Coal&1800&186&146.22&0.347\\
Oil&400&41&211.64&0.29\\
CCGT&900&94&85.60&0.2\\
Hydro&2500&255&0.00&0\\
Wind&1137&102&0.00&0\\
Solar&395&34&0.00&0\\
\hline
\end{tabular}%
\end{table}

For capacity market derating factors ($DR_{t}$), we assume a value of zero for wind and solar, and a value of one for all other technologies. This is because we assume hydro and fossil-fuel based technologies are always available for generation while wind and solar generation are, at times, not available at all for generation. 

\subsubsection{Renewable power generation data}
\label{sec:input-data}
As with \cite{lynch2019role} and \cite{devine2023role}, data from the MERRA2 reanalysis \citep{Bosilovich2016} were used to generate input wind and PV data. Note that wind and solar PV are not only variable but also uncertain and their uncertainties are correlated since both depend on the meteorological conditions. It is therefore important to take these correlations into account when providing input data for the stochastic MCP. 

To calculate, wind capacity factors ($NORM^{\text{WIND}}_{t,p,s}$), the method described in \cite{Cradden2017} was used to convert MERRA2 wind data to wind power. Similarly, to calculate hourly PV capacity factors ($NORM^{\text{PV}}_{p,s}$), the MERRA2 variables SWGDN and T2M were converted into PV power using a phyiscal PV model \citep{ritzenhoff1992erstellung, ruppert2016impact, schwarz2018self, schwarz2018two}\footnote{All values for $NORM^{\text{WIND}}_{t,p,s}$ and $NORM^{\text{PV}}_{p,s}$ are provided in the online appendix\footnotemark[1]. The values for $NORM^{\text{WIND}}_{t,p,s}$ are split into three different wind regions in Ireland - see \cite{bertsch2018analysing} for details.}. The hourly wind and solar capacity factor time series of 35 years were then clustered into six representative years. Thus we consider $|S|=6$ scenarios for the MCP model described in Section \ref{sec:methodology}. The renewable data generation and the clustering procedure are described in further detail in \cite{bertsch2018analysing} while the probabilities and chosen years of occurrence are summarised in Table \ref{tab:scen_prob}. To ensure that the spatial and temporal correlations between wind and PV are preserved, we use these historical wind and PV data as a basis for our analysis.  We consider 250MW of solar PV capacity, 125MW split amongst for the generating firms (Table \ref{tab:firm_data}) and 125MW split evenly between amongst the two prosumer consumer groups. Furthermore, Table \ref{tab:firm_data} shows that we consider an initial 5600MW of installed wind capacity. We base these assumptions around installed RES-E capacity estimates for 2030 from \citet{eirgrid2019}.

\begin{table}[htp!]
 \centering
\footnotesize
\caption{Representative years chosen for RES scenarios and corresponding probabilities of occurrence; see \cite{bertsch2018analysing}}\label{tab:scen_prob}
\begin{tabular}{lrrrrrr}
\hline
Year & 1983 & 1998 & 2001 & 2003 & 2004 & 2015 \\
Probability of occurrence ($PROB_{s}$) & 0.486 & 0.286 & 0.086 & 0.086 & 0.029 & 0.029 \\
\hline
\end{tabular}%
\end{table}

We assume a FiP payment of $X_{t}$ = 23 \euro/MWh for wind and PV generation. These values are obtained from \cite{farrell2017specifying}. A FiP is not considered  for any other technology.

\section{Results} \label{sec:results}
In this section, we present the results of our study. In Section \ref{sec:results_MP} we study the interaction between market power and investment decisions on the demand  side and the supply side. As a sensitivity analysis, we analyse the effects of the absence vs. presence of revenues from a Feed-in Premium for renewables in Section \ref{sec:results_FiP}. In Section \ref{sec:results_CO2} we examine carbon emissions for each of the cases presented. In cases where market power is present in the model, all generating firms are price-makers, modelled \emph{\'{a} la} Cournot, while in cases where no market power is present, all firms are modelled as price-takers, i.e., the market is perfectly competitive. Throughout all analyses, we consider four different levels for the proportion of prosumers in the model, 0\%, 33\%, 67\%, and 100\%.

\subsection{Interaction of market power and investment decisions}\label{sec:results_MP}
Throughout Section \ref{sec:results_MP}, we assume that the generators receive a Feed-in Premium for renewable generation in addition to the energy market revenues.

\subsubsection{Demand side effects}

\begin{figure}[h]
     \centering
     \begin{subfigure}[b]{0.475\textwidth}
         \centering
       \epsfig{figure=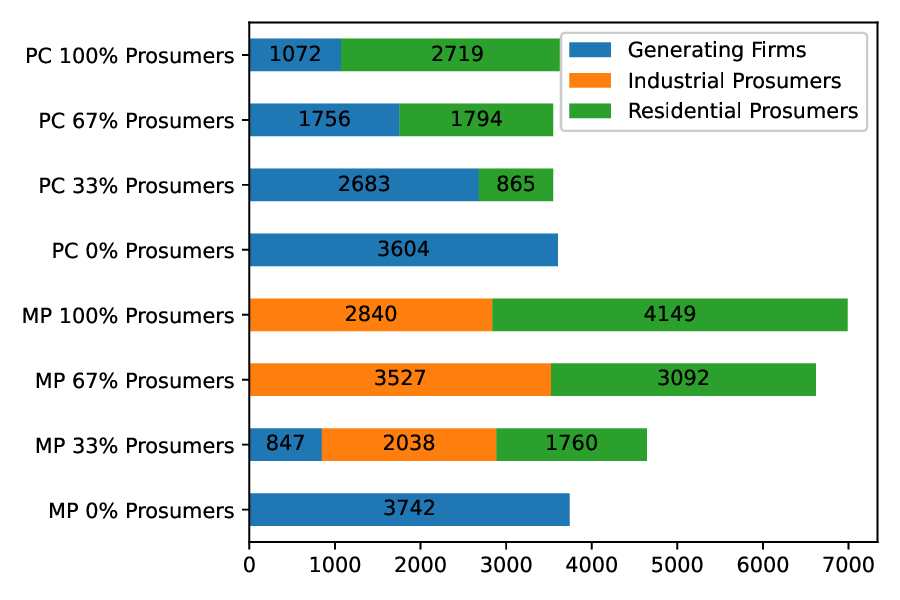,width=\textwidth}
         \caption{PV investments (MW)}
         \label{fig:invest_pv}
     \end{subfigure}
     \hfill
     \begin{subfigure}[b]{0.475\textwidth}
         \centering
        \epsfig{figure=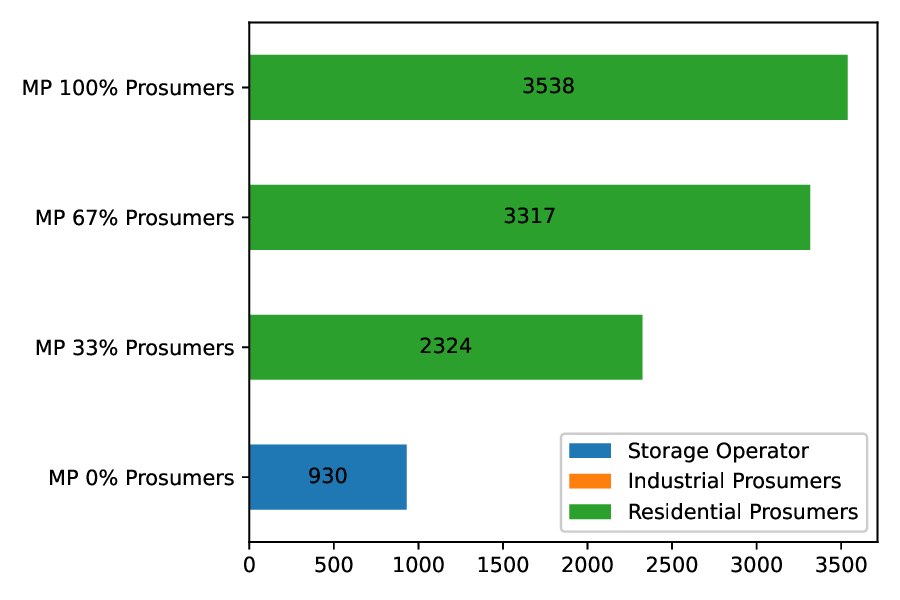,width=\textwidth}
         \caption{Storage investment (MW); market power cases only}
         \label{fig:invest_storage}
     \end{subfigure}
        \caption{Self-sufficiency investments under Market Power (MP) and Perfect Competition (PC)}
        \label{fig:invest_impact_of_MP_demand}
\end{figure}

We start by examining the interaction between market power and self-sufficiency investments. In addition, we also discuss the effects market power has on consumer tariffs. Figures \ref{fig:invest_pv} and \ref{fig:invest_storage} display PV and storage investments by residential and industrial consumers, respectively. 

In the absence of market power, Figure \ref{fig:invest_pv} shows that residential prosumers invest in PV with the number of investments increasing as the \% of prosumers increases, as one could expect. However, industrial prosumers do not invest in PV in the absence of market power. This is mainly driven by the higher ``retail premium'' (difference between wholesale and retail prices) of residential prosumers as compared to industrial consumers, which basically makes PV-based self-sufficiency more attractive for residential consumers despite the slightly higher technology costs of residential-scale PV investments. 

There are zero investments in storage when market power is absent. While more than 2GW of conventional generation capacity are expected to be decommissioned in both cases (with and without market power), a significant amount of flexible and dispatchable conventional generation capacity will remain available in the system. Therefore, in the absence of market power, there is no need in the studied system for additional flexibility in the form of storage. And since prices in the absence of market power are comparatively low (see Figure \ref{fig:prices_impact_FiP} on page \pageref{fig:prices_impact_FiP}), there is also no sufficient incentive for prosumers to invest in storage. 

In the presence of market power, Figure \ref{fig:invest_pv} shows substantially increased investments in PV by residential prosumers. In addition, in contrast to the perfect competition case, industrial prosumers also now invest in PV. As the \% of prosumers goes from 33\% to 67\%, PV investments by industrial prosumers increase from 2038MW to 3527MW. However once all consumers are prosumers, industrial prosumers reduce their PV investments to 2840MW, despite extra industrial consumers with the opportunity to invest in self-sufficiency. This is because there are also further residential prosumers who increase the level of PV investments in that group, which leads to reduced system prices and consequently reduced profitability of PV-based self-sufficiency for industrial prosumers.

Figure \ref{fig:invest_storage} displays investments in battery storage for the market power cases. When there are no prosumers, only the storage operator invests into storage (in addition to the initial 700MW it begins with). However, once prosumers enter the market, residential prosumers begin to invest in battery storage. Because all residential consumers buy electricity at the retail price, they have more savings to gain from increased self-sufficiency by investing in storage compared to the independent storage operator who generates margins from arbitrage in the wholesale market -- particularly in the presence of market power. These prosumers use storage to shave peak prices and also to increase their self-sufficiency, which helps them reduce the volume of and costs for electricity consumed from the grid. This reduces profitability for the independent storage operator who consequently does not invest in further storage. While only residential prosumers invest in storage, all consumers and prosumers still benefit from the reduced prices, particularly at peak times.

\begin{figure}[h]
     \centering
     \begin{subfigure}[b]{0.45\textwidth}
         \centering
         \epsfig{figure=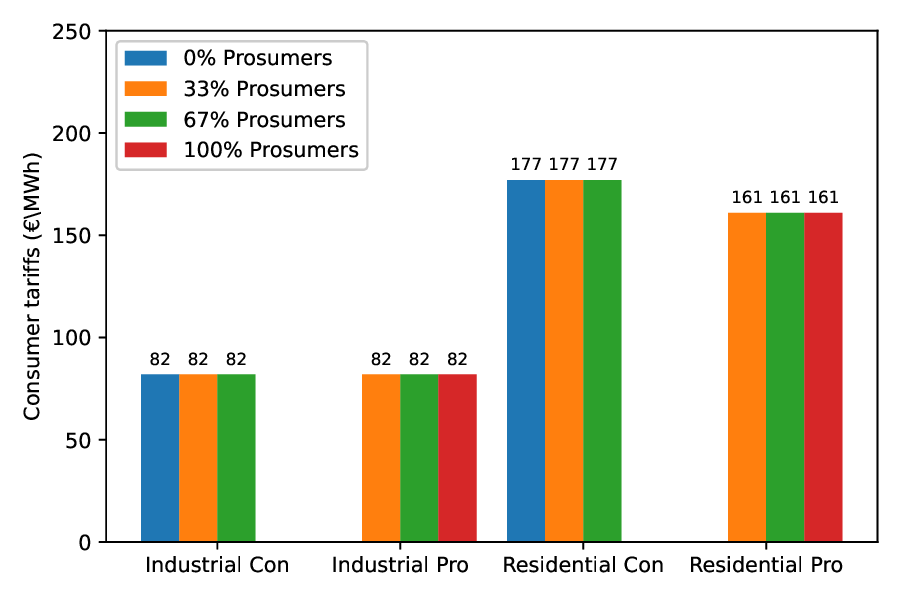,width=\textwidth}
        \caption{Perfect competition}
         \label{fig:tariffs_cm}
     \end{subfigure}
     \hfill
     \begin{subfigure}[b]{0.45\textwidth}
         \centering
         \epsfig{figure=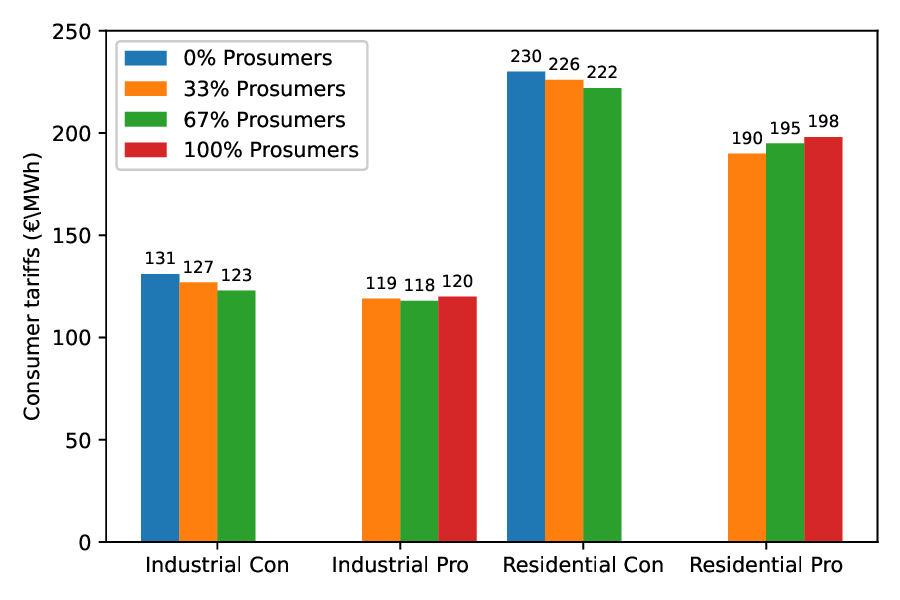,width=\textwidth}
         \caption{Market Power}
         \label{fig:tariffs_mp}
     \end{subfigure}
     \caption{Consumer tariffs (\euro/MWh ) for traditional consumers (Con) and prosumers (Pro)}\label{fig:tariffs_full}
\end{figure}

Figure \ref{fig:tariffs_full} displays the resulting consumer tariffs. They are calculated for each consumer group by dividing their objective function value (equation \eqref{eqn:obj_act}) by their reference demand ($D^{\text{REF}}_{k,p}$). In perfect competition, there is no difference between the traditional industrial consumers and the industrial prosumers. This is because the industrial prosumers do not invest in self-sufficiency, despite having the opportunity to do so. In contrast, residential prosumers do invest in PV under perfect competition and, consequently, their tariffs reduce as compared to the residential consumers who do not prosume. 

When market power is assumed, tariffs increase for all consumers, as expected. For industrial consumers, the ability to prosume reduces their tariffs by  3-13 \euro/MWh, depending on the number of prosumers assumed. For residential consumers, savings offered by prosuming are 24-40 \euro/MWh. These savings are greater than the reductions in wholesale prices observed in Figure \ref{fig:prices_impact_FiP}. This is because of investment in battery storage by the prosumers, enabling them to move their demand to lower priced periods of time. In addition, self-sufficiency allows the prosumers to avoid paying the retail price including the premium on top of the wholesale price. 

Interestingly, as the \% of prosumers increases, the tariff for traditional consumers decreases. This shows that those who do not invest in self-sufficiency still benefit from the investment of other consumer groups. That is, investments by prosumers reduce the price-making ability of the firms and hence lead to reduced prices for all. 

Moreover, as the \% of prosumers increases, the tariff for prosumers increases, particularly in case of residential prosumers. This suggests that when the \% of prosumers are low, prosumers can make bigger savings but these specific savings in \euro/MWh somewhat diminish as further consumers become prosumers. In other words, while investments in self-sufficiency increase as the \% of prosumers increases, the marginal benefits of these additional investments decrease. 

Finally, we also study the impact of market power on load shedding. We find that the increased prices under market power make load shedding an economically attractive option and that the ability to prosume can reduce total annual load shedding by around 8\%. However, we also find that the ability to prosume does not affect the max hourly load shedding.

\subsubsection{Supply side effects}

\begin{figure}
      \centering
    \epsfig{figure=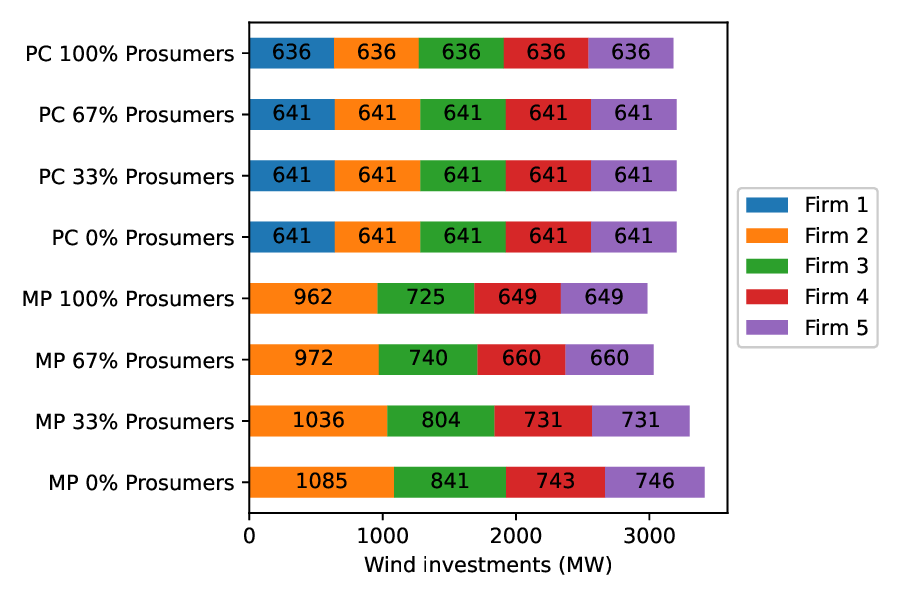,width=0.45\textwidth}
      \caption{Wind investments under Market Power (MP) and Perfect Competition (PC)}
    \label{fig:invest_wind}
\end{figure}

In this section we examine the impacts market power and self-sufficiency investments have on the generation firms’ investment decisions, their profits, and wholesale prices.  Figure \ref{fig:invest_wind} displays investments for wind in greater detail by firm. For each of the cases considered, we do not observe any investments in any other generating technology (Coal, Oil, CCGT, Hydro) despite the firms having the ability to do so. 

In the absence of market power, total wind investment is similar across the different \% of prosumers cases considered. For each test case, the firms invest the same amount as each other in perfect competition as wind investment costs are the same for each firm. 

While PV investments by prosumers increase as the \% of prosumers increases, total PV investments are also relatively constant across the \% of prosumers test cases, in the absence of market power. In the case where there are no prosumers, only the generating firms invest in PV of course. As the \% of prosumers increases, the firms decrease their investment in PV while residential consumers increase theirs.

\begin{figure}[h]
     \centering
     \begin{subfigure}[b]{0.45\textwidth}
         \centering
          \epsfig{figure=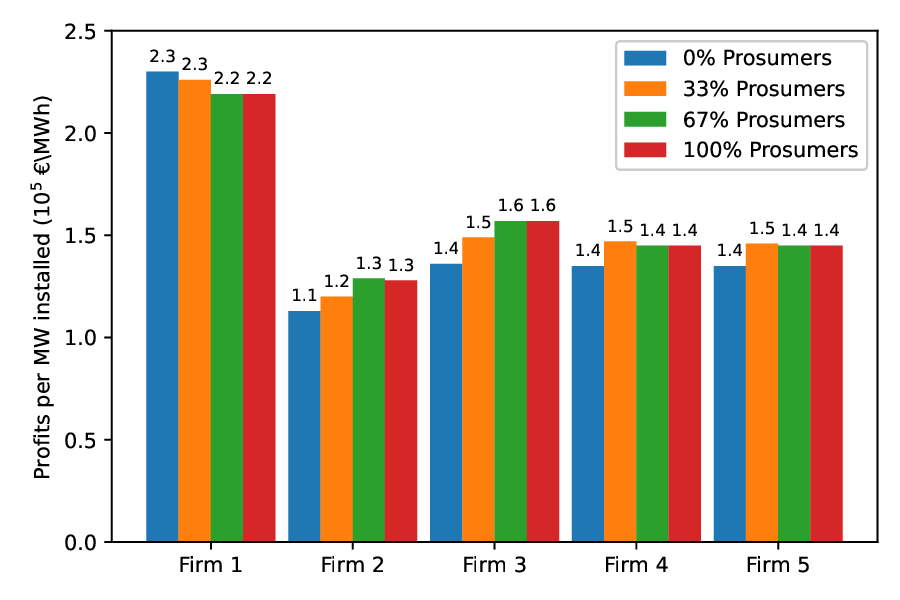,width=\textwidth}
         \caption{Market Power}
         \label{fig:profits_mp_mw}
     \end{subfigure}
     \hfill
     \begin{subfigure}[b]{0.45\textwidth}
         \centering
          \epsfig{figure=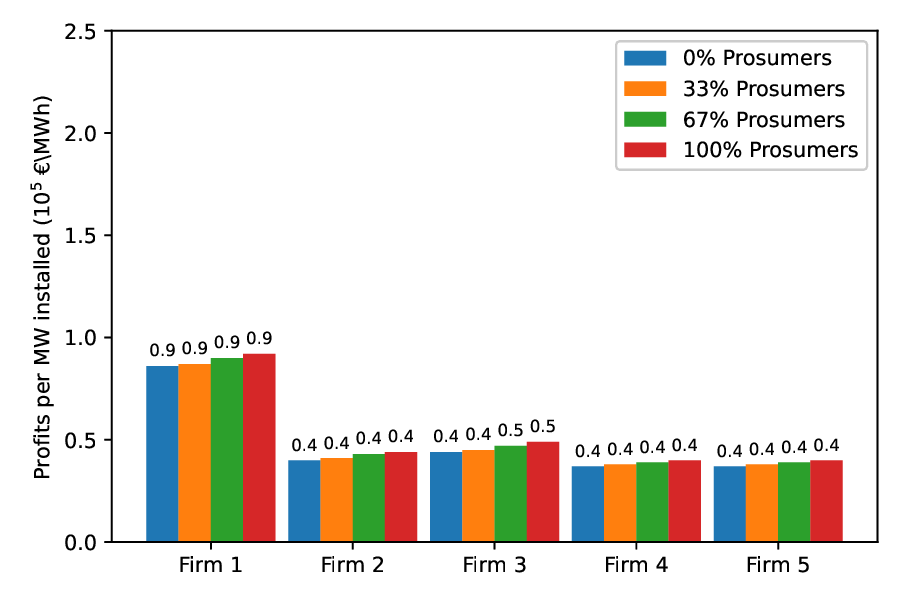,width=\textwidth}
         \caption{Perfect competition}
         \label{fig:profits_cm_mw}
     \end{subfigure}
     \caption{Firm profits per MW installed (\euro hundreds of thousands) }\label{fig:profits_impact_mp_mw}
\end{figure}

In the presence of market power, Figures Figure \ref{fig:invest_pv} and \ref{fig:invest_wind} shows substantially increased investments in PV but not for wind. More specifically, when there are no prosumers, the presence of market power slightly increases total wind investment (3203MW to 3415MW). However, the increased PV investments under market power resulting from the increased \% of prosumers, lead to decreased investments in wind. So much so that when 100\% of consumers are prosumers, there are less wind investments (2985MW) compared to the equivalent perfect competition case (3180MW). The decreased prices resulting from increased PV lead to less investments in wind by the generation firms. Although, it is important to note that, despite reduced investments by the firms, their profits - as expected - are substantially higher when market power is present (Figure \ref{fig:profits_impact_mp_mw}). 

In contrast to the perfect competition cases, there is heterogeneity in the firms' wind investment. In the presence of MP, firm $f=1$ does not invest in wind. In turn firm $f=2$'s investment in wind are the highest, followed by firm $f=3$. Firm $f=4$ and firm $f=5$'s investment in wind are similar. The firms' investment in wind is inversely proportional to their initial wind generation capacity (Table \ref{tab:firm_data}). This is because we model the firms' market power \emph{\'{a} la} Cournot, where firms reduce their generation levels in order to increase the system price. 

Firm $f=1$ has the largest generating portfolio and thus strategically reduces its generation the most. Consequently, if firm $f=1$ was to invest in wind, it would not use that generation all that often. Hence, in the presence of market power, it is not profitable for firm $f=1$ to invest in wind. Similarly, the next largest firms ($f=4$ and $f=5$) do not invest as much in wind as firms $f=2$ and $f=3$ as they do would not use it as much, due to similar strategic behaviour. Another effect of assuming all generating firms behave  \emph{\'{a} la} Cournot competition is that we observe an equalizing effect. That is, the difference in the size of the firms' portfolios post-investment is smaller than pre-investment. 

Figure \ref{fig:profits_impact_mp_mw} displays the firms' profits per MW installed. Firm $f=1$ has the largest profit per MW installed. This is because firm $f=1$ has the largest portfolio of renewable generators (wind, PV, and hydro; see Table \ref{tab:firm_data}), which make the largest profits due to their low variable power generation costs. 

In the absence of market power, each firms' profits per MW installed increase as the \% of prosumers grows. This is because firms decrease investment in PV as more prosumers enter the market. While this has little effect on the firms' absolute profits, it increases their profits per MW installed. 

As expected, once market power is assumed in the model, profits for all firms increase substantially. Again, once prosumers enter the market, profits of firms 2 - 5 increase. However, this increase is not monotonic as these firms' capacity after investment and decommissioning changes very little, as the \% of prosumers grows from 33\% to 100\%.

The exception to this is firm $f=1$'s profits under market power. Firm $f=1$ does not invest in PV when market power is present, regardless of the number of prosumers. Hence, as the \% of prosumers increases, firm $f=1$'s profits decrease due to declining prices.

The absolute profits of all generating firms increase when the assumption of perfect competition changes to market power, regardless of what (fixed) percentage of prosumers is assumed.  When the level of prosumers increases, the profit of all firms (absolute) decreases, regardless of whether market power is assumed or not.

\begin{figure}[h]
    \centering
    \epsfig{figure=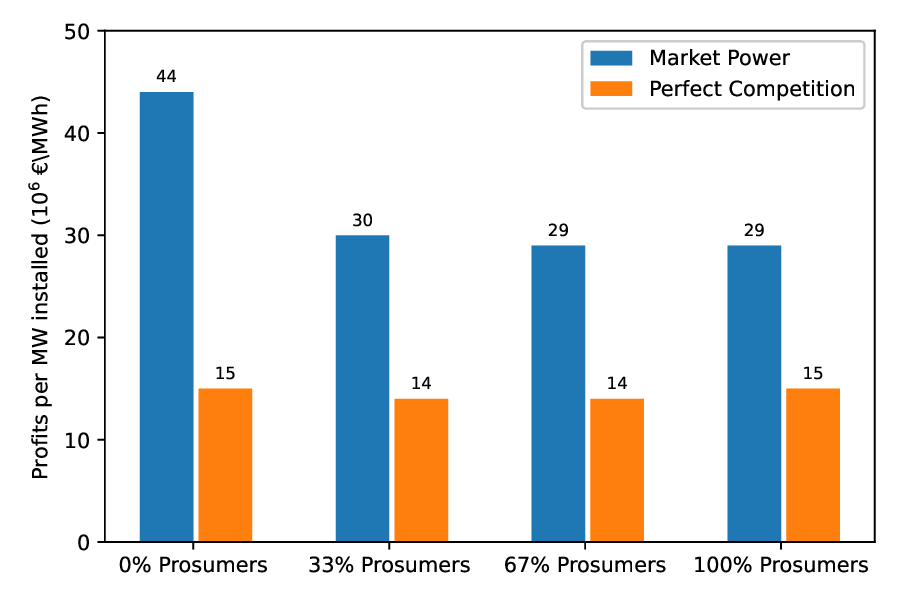,width=0.45\textwidth}
    \caption{Storage operator's profits (\euro Millions)}
    \label{fig:impact_mp_stor_profits}
\end{figure}

Figure \ref{fig:impact_mp_stor_profits} displays the independent storage operator's profits. In perfect competition, they do not invest into further storage. The profits they make are from the initial 700MW of battery storage we assume they have (Section \ref{sec:data_SO}). Moreover, in perfect competition, the \% of prosumers has little effect on the storage operator's profits. This is because the storage operator generates margins from arbitrage in the wholesale market and, under perfect competition, the presence of prosumers does not have a massive impact on wholesale prices. 

When market power is present, the storage operator's profits increase by \euro 14 to \euro 29 million, depending on the number of prosumers. However, the arrival of prosumers reduces the storage operator's profits. On the one hand, the self-sufficiency by prosumers leads to reduced price peaks, which explains reduced storage profits. On the other hand, as shown in Figure \ref{fig:invest_storage}, investment in storage switches from the storage operator to residential prosumers when market power is present.

With regard to exit decisions, we observe that all coal generation is decommissioned in all test cases considered. The presence of market power increases the decommissioning of CCGT generation, particularly for firm $f=1$, which consequently leads to decreased decommissioning of oil-fired generation. This is because when firm $f=1$ behaves strategically, it utilises its CCGT generation less often and because oil-fired generation has cheaper maintenance costs than CCGT.

\subsection{Sensitivity analysis: removal of FiP}\label{sec:results_FiP}
In recent years, the capital costs associated with renewables have been decreasing, thus questioning the need for Feed-in Premia \citep{finke2023exploring}.
In this section, we examine the effect that the FiP has on investments, prices, consumer tariffs, and load shedding. 

\begin{figure}[h]
     \centering
     \begin{subfigure}[b]{0.475\textwidth}
         \centering
         \epsfig{figure=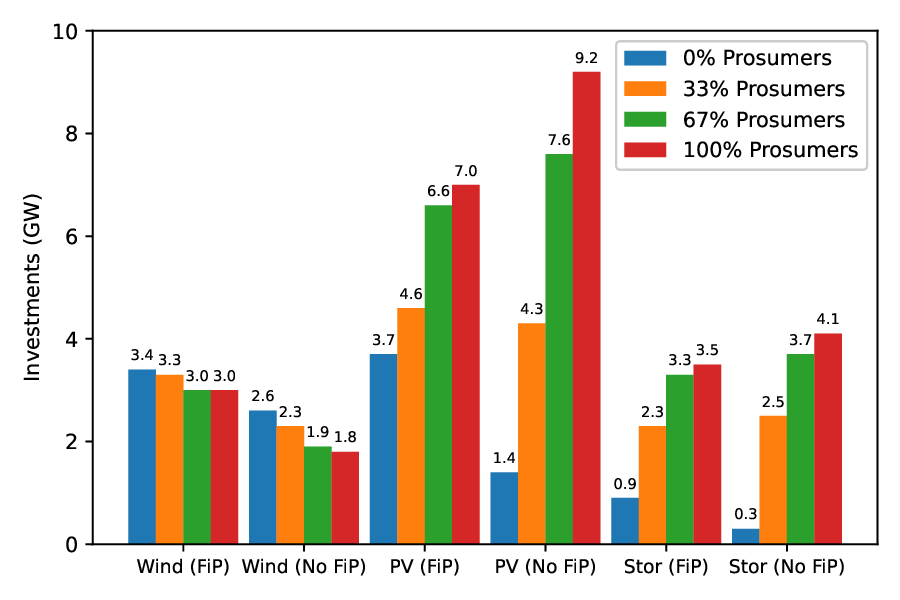,width=\textwidth}
         \caption{Market power}
         \label{fig:invest_mp_impact_fip}
     \end{subfigure}
     \hfill
     \begin{subfigure}[b]{0.475\textwidth}
         \centering
       \epsfig{figure=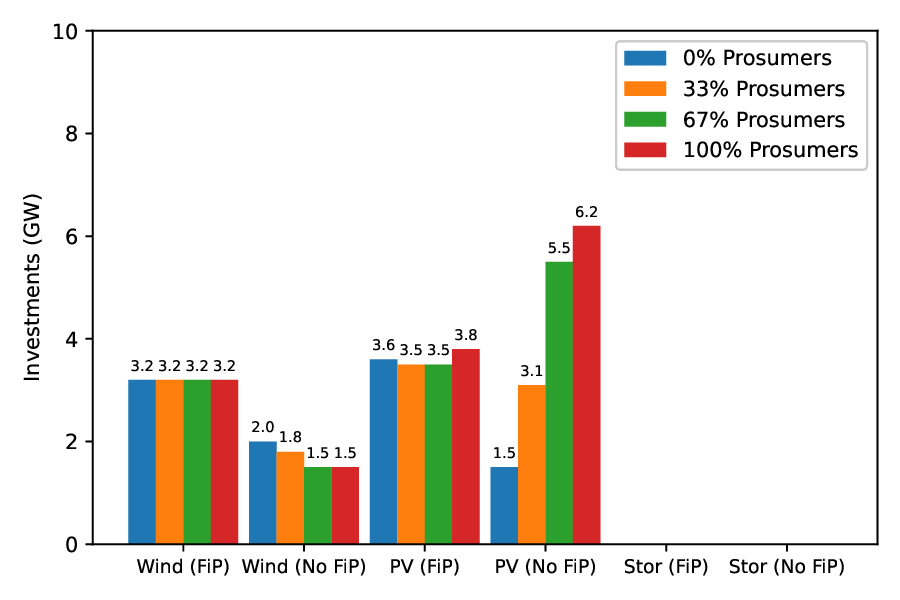,width=\textwidth}
         \caption{Perfect competition}
         \label{fig:invest_cm_impact_fip}
     \end{subfigure}
             \caption{Impact of FiP on system-wide investments in storage (stor), wind, and PV (GW)}\label{fig:invest_impact_fip}
\end{figure}

Figure \ref{fig:invest_impact_fip} shows, as one might suspect, that the presence of a FiP increases wind investment. And this holds regardless of the number of prosumers. In contrast, the presence of a FiP has a mixed effect on PV investments. When there are no prosumers, the FiP increases PV investments from the generating firms. The FiP increases the firms' revenues and hence the profitability associated with such investments. 

\begin{figure}[h]
    \centering
    \epsfig{figure=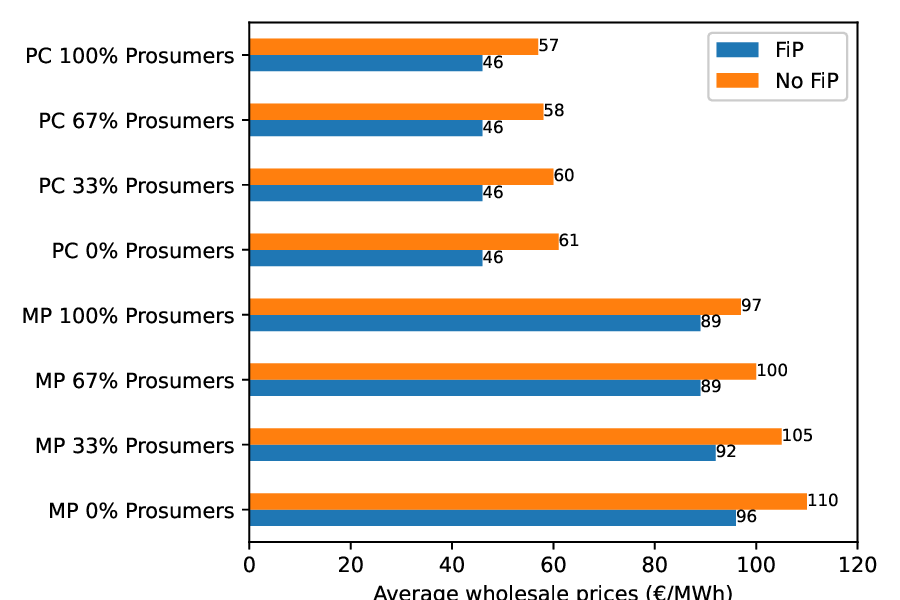,width=0.5\textwidth}
    \caption{Average wholesale prices (\euro /MWh) under Market Power (MP) and Perfect Competition (PC)}
    \label{fig:prices_impact_FiP}
\end{figure}

However, the arrival of prosumers means PV investment starts to switch from the generating firms to the prosumers. When there is a majority of prosumers (67\% and 100\%), the absence of the FiP increases PV investment compared to a case where a FiP is available. This result appears counter-intuitive. However, the absence of a FiP decreases wind investments and hence increases wholesale prices (Figure \ref{fig:prices_impact_FiP}). Consequently, the increased prices lead to increased profitability of PV-based self-sufficiency with the level of PV investments increasing as the \% of prosumers grows. As Figure \ref{fig:invest_impact_fip} shows, this result holds regardless of whether the market is perfectly competitive or not. 

The effect that the FiP has on storage investments is correlated with the effect on PV. That is, as PV investment increases, so does storage. When prosumers are present in the market, both investments are mainly observed on the demand side, where storage investment can only be observed when market prices are elevated due to market power. Particularly under market power, such investments enable prosumers to cover a relatively large share of their electricity demand by own generation, which helps protect them against high market prices.

Figure \ref{fig:prices_impact_FiP} displays average wholesale prices with and without the FiP. The FiP decreases wholesale prices by \euro 9 - \euro 15 /MWh in perfect competition and by \euro 8 - \euro 14 /MWh when market power is present. This is due to increased wind investments. The effect is largest when there are no prosumers - as discussed above. In these cases, the FiP also increases PV investments. As the \% of prosumers grows, the size of the effect decreases. This is because of the decreases in PV investments from prosumers when the FiP is present (see Figure \ref{fig:invest_impact_fip}), which can be explained by the slightly reduced attractiveness of self-sufficiency investments resulting from higher wind investments which in turn lead to reduced wholesale prices.

\begin{figure}[h]
     \centering
     \begin{subfigure}[b]{0.475\textwidth}
         \centering
         \epsfig{figure=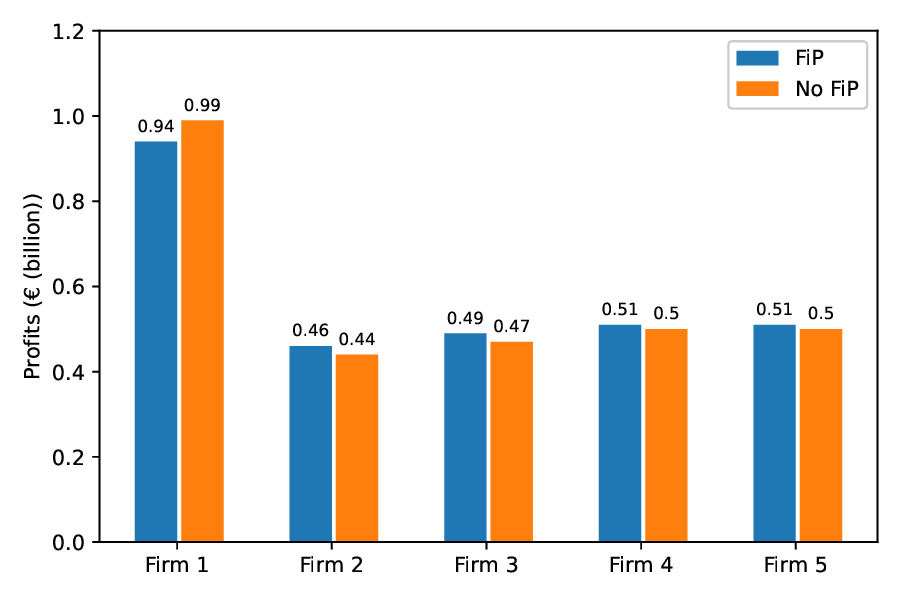,width=\textwidth}
         \caption{0\% Prosumers}
         \label{fig:iprofits_fip_0P}
     \end{subfigure}
     \hfill
     \begin{subfigure}[b]{0.475\textwidth}
         \centering
       \epsfig{figure=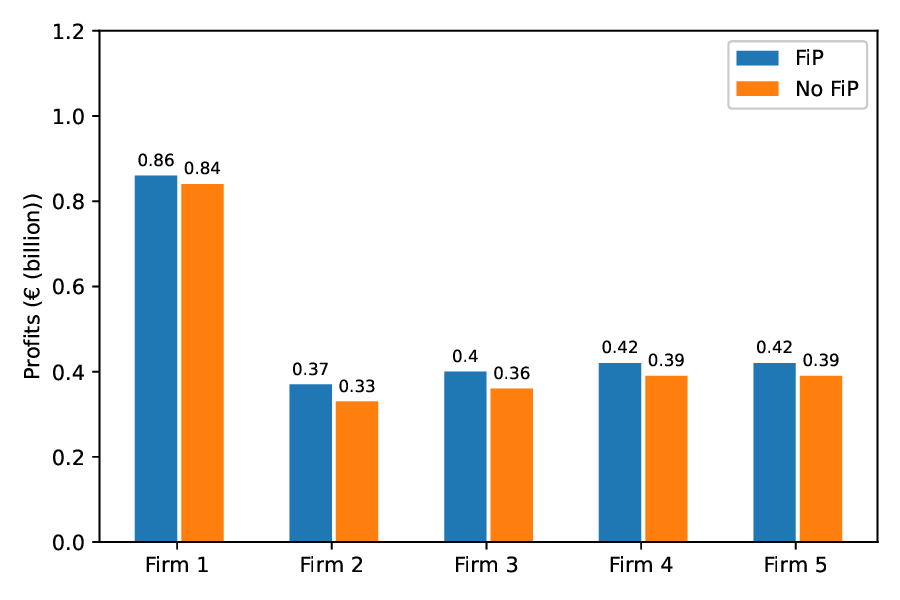,width=\textwidth}
         \caption{100\% Prosumers}
         \label{fig:iprofits_fip_100P}
     \end{subfigure}
\caption{Impact of FiP on firms' profits (\euro billions); market power cases only} \label{fig:profits_impact_fip}
\end{figure}

Figure \ref{fig:profits_impact_fip} examines the effect that the FiP has on the firms' profits. While the absolute profits are high, the relative changes resulting from the FiP are generally rather small. In both the absence and presence of market power, the FiP increases profits for firms $f=2$ - $f=5$ as the FiP provides extra revenues to these firms. 

The FiP has a mixed effect on the largest firm, firm $f=1$. When there are no prosumers, the FiP decreases firm $f=1$'s profits. This is because the FiP encourages the other firms to invest in wind generation and hence this reduces firm $f=1$'s market share. Firm $f=1$ does not invest in wind when the firms are assumed to have market power. 

However, as the \% of prosumers increases, the impact of the FiP on firm $f=1$'s profits flips. PV investments by prosumers increase as the \% of prosumers grows. The presence of a FiP switches investments back from prosumers to firms $f=2$ - $f=5$; increased wind from the firms and decreased PV investment from the prosumers (Figure \ref{fig:invest_impact_fip}). In contrast to the prosumers, the firms behave \emph{\'{a} la} Cournot with this generation. This means that the relative decrease in wholesale prices due the FiP is not as big when there is a majority of prosumers. This benefits firm $f=1$. 

When there are no prosumers, the absence of a FiP increases the generating firms' profit per MW installed. This is because the firms invest in less capacity when there is no FiP. Increasing the \% of prosumers leads to the absence/presence of a FiP having little effect on the firms' profit per MW installed.

Finally, we analyse the effect that the FiP has on consumer tariffs. As expected, we find that the FiP reduces consumer tariffs across all levels of prosumers and regardless of whether market power is present or not. Note, however, that this finding does not account for who pays for the FiP. That is, is it electricity consumers via their bills or is it paid via general taxation through the exchequer. In addition, it is important to note that the impact of the FiP on consumer tariffs is much smaller than that of market power.

\subsection{Impact of market power, self-sufficiency and FiP on Carbon Emissions}\label{sec:results_CO2}
\begin{figure}[h]
    \centering
    \epsfig{figure=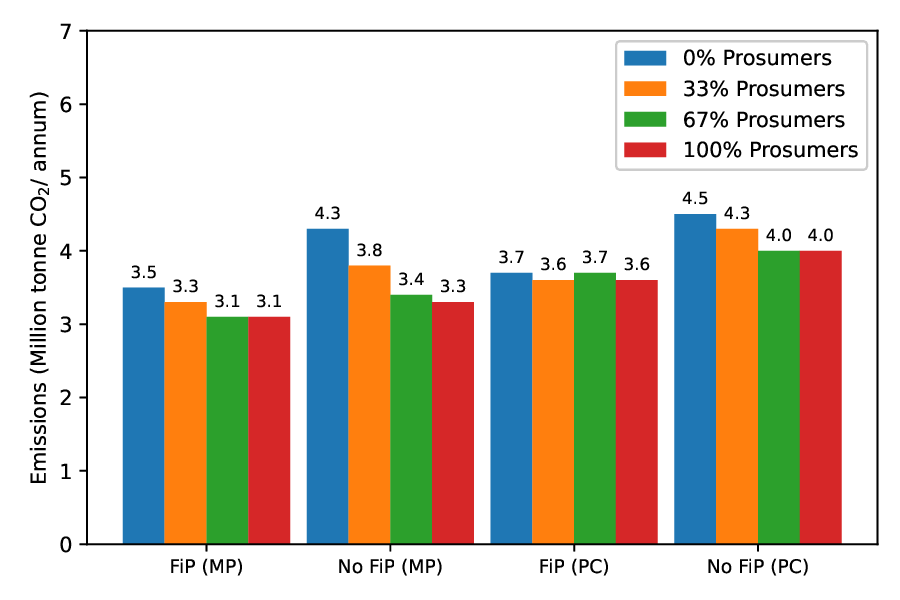,width=0.475\textwidth}
    \caption{Carbon emissions (Million tonne CO$_{2}$/annum) under Market Power (MP) and Perfect Competition (PC) in presence (FiP) and absence (No FiP) of a FiP}
    \label{fig:C02}
\end{figure}

In this section, we examine the effects that both market power and the Feed-in Premium have on carbon emissions. Figure \ref{fig:C02} shows that market power decreases carbon emissions by 6\% - 22\%, depending on the test case considered. This effect is present when the FiP is included and excluded from the model. Market power increases investment in PV when the FiP is present (Figure \ref{fig:invest_pv}) and increases investment in both PV and wind when the FiP is not present (Figure \ref{fig:invest_mp_impact_fip}). 

In the presence of market power, as the \% of prosumers and hence PV and storage investments increase, emissions decrease. PV and storage enable prosumers to cover a relatively large share of their demand by own, climate-friendly electricity generation, which in the end displaces fossil-fuel generation on the supply side.

Similarly, without market power and without the FiP, carbon emissions decrease as the \% of prosumers increases. However, without market power but with FiP, the \% of prosumers does not affect carbon emissions. With the FiP, there are high levels of investment in wind (3200MW) which lead to relatively low prices (Figure \ref{fig:prices_impact_FiP}). Without the FiP, investment in wind decreases by ($\Tilde{}$ 1800MW). This leads to higher prices and consequently increased levels of PV investments for self-sufficiency, which increases as the \% of prosumers increases, resulting in lower emissions (Figure \ref{fig:invest_cm_impact_fip}).  

\begin{figure}[h]
     \centering
         \centering
       \epsfig{figure=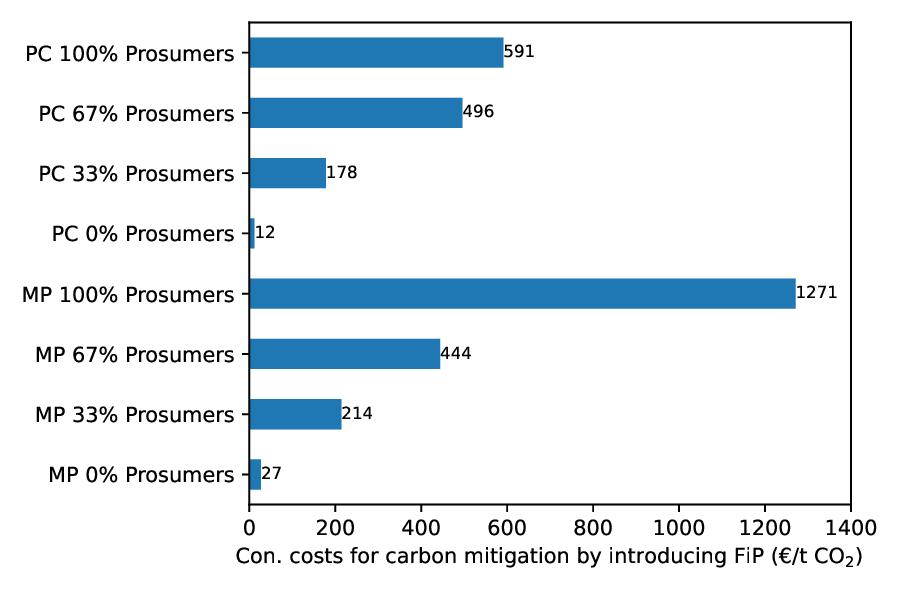,width=0.7\textwidth}
     \caption{Consumer costs for carbon mitigation by introducing a FiP in \euro/t CO$_2$ under market power (MP) as well as perfect competition (PC)}\label{fig:CarbonMitigationConsumerCostByFiP}
\end{figure}

Figure \ref{fig:CarbonMitigationConsumerCostByFiP} puts the total additional consumer costs from introducing a FiP into perspective with the emissions saved (increase in consumer costs divided by decrease in emissions). As opposed to our analysis of consumer tariffs in Section \ref{sec:results_FiP}, the total consumer costs here include the capacity costs and FiP costs borne by the consumers in addition to the consumer objective. The results in Figure \ref{fig:CarbonMitigationConsumerCostByFiP} reveal that, in the absence of any prosumers, the introduction of a FiP enables carbon emission reduction at rather low costs to consumers (27\euro/t CO$_2$ in the presence and 12\euro/t CO$_2$ in the absence of market power).\footnote{For comparison: Actual CO$_2$ market prices in the EU's ETS varied roughly within a range of 50-100 \euro/t CO$_2$ in 2023 and the beginning of 2024.} As soon as consumers are able to invest in PV and storage, however, the additional benefits in terms of emission reduction that can be achieved by introducing a FiP become increasingly expensive. This effect is even more pronounced when market power is present. This is because the demand side investments in PV and storage lead to carbon emission reductions and, as discussed above, in the absence of a FiP and in the presence of market power, market prices are highest, which makes these self-sufficiency investments economically very attractive for consumers. As described above, the introduction of a FiP mainly leads to higher investments in wind power, which displaces PV investments on the demand side. That is, when demand side investments are possible, the introduction of a FiP leads to a displacement of one low-carbon technology by another one.

\subsection{Quantifying further effects of self-sufficiency investments: a back of the envelope analysis}
Besides the interactions between market power, FiP, self-sufficiency investments and emissions presented in the previous subsections, there are further implications involved. Naturally, self-sufficiency leads to reduced grid demand by prosumers. In most jurisdictions, all consumers and prosumers do not just pay the wholesale price for their electricity demand from the grid, but typically pay different (usually volumetric) charges (e.g., grid charges, as well as other levies, taxes, or procurement margins of retailers). In this paper, we have considered these charges in a simplified way as a volumetric lumpsum of 35 \euro/MWh for industrial and 125 \euro/MWh for residential consumers, which we called `retail premium'. It is important to note that these charges constitute a crucial element for recovering the costs for maintaining the electricity system as a whole. Without changing the way in which these charges are collected (i.e. when they remain volumetric charges), this means that reduced grid demand resulting from self-sufficiency also leads to a cost recovery gap for the system infrastructure \citep[cf. e.g. ][]{bertsch2017drives,schwarz2018self,mehigan2018review,aniello2023shaping}. 

\begin{figure}[h]
    \centering
    \epsfig{figure=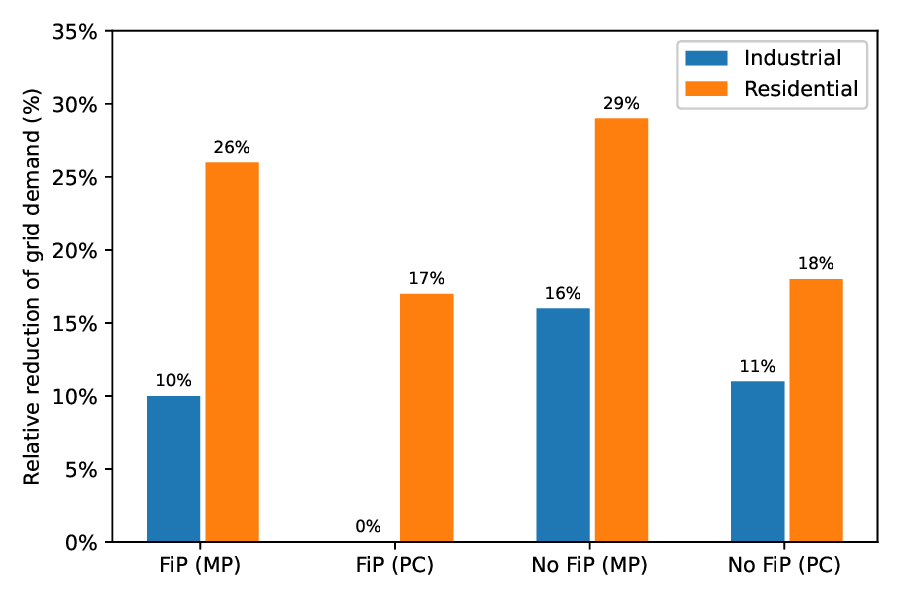,width=0.6\textwidth}
    \caption{Relative reduction of grid demand following from self-sufficiency under cases with market power (MP) and perfect competition (PC), each in presence and absence of a Feed-in Premium (FiP)}
    \label{fig:BOE}
\end{figure}

For cases, where 100\% of consumers are able to invest in PV and/or storage, Figure \ref{fig:BOE} shows the relative grid demand reduction that can be expected. Relative grid demand reduction is defined as follows:
\begin{subequations}
\begin{equation}
    \frac{\sum_{k,p}\big(D^{\text{REF}}_{k,p}-D^{\text{GRID}}_{k,p}\big)}{\sum_{k,p}D^{\text{REF}}_{k,p}},
\end{equation}
where consumer group $k$'s grid demand is their reference demand less any self-consumption:
\begin{equation}
 D^{\text{GRID}}_{k,p}=D^{\text{REF}}_{k,p}+g^{\text{UP}}_{k,p,s}-(1-LOSS_{k})g^{\text{DOWN}}_{k,p,s}-g^{\text{PV}}_{k,p,s}.
\end{equation}
\end{subequations}
In line with the findings presented so far, self-sufficiency and hence grid demand reduction are higher when market power is present than under perfect competition and they are higher for residential compared to industrial consumers. Overall the grid demand reduction amounts to roughly 15-30\% for residential consumers and to roughly 10-15\% for industrial consumers. When multiplying the absolute grid demand reduction by the retail premia for residential and industrial consumers respectively, this results in a total cost recovery gap (industrial and residential consumers together) of roughly 290-620 Million \euro/a depending on whether market power and a FiP are present or not. As a reference point, in the cases, where no prosumers are present in the system at all, i.e. the retail premia would need to be paid for the complete reference demand, the entire retail premia collected would amount to roughly 2.5 Billion \euro/a. 

\section{Discussion, conclusions and outlook} \label{sec:discussion}

In this work, we introduce a stochastic Mixed Complementarity Problem to model an electricity market characterised by an oligopoly. On the supply side, we consider a number of generating firms, each of whom hold and may invest in multiple generating technologies. On the demand side, we consider both industrial and residential consumers. Furthermore, we distinguish between traditional consumers and consumers who may prosume, that is prosumers you can invest in self-sufficiency such as PV and battery storage. In addition, we consider an independent storage operator.

We go beyond the existing literature by studying the interactions between market power, self-sufficiency investments on the demand side, and supply-side investments. That is, we analyse how market power affects self-sufficiency investments and, in turn, if self-sufficiency can help mitigate against market power. Moreover, we carry out our analyses in a multi-revenue stream setting where generators receive remuneration in both an energy and capacity market in addition to receiving a FiP for renewables. 

We study how both market power and self-sufficiency impact investment decisions made by generation firms. Furthermore, we study how self-sufficiency impacts each of the considered consumer groups differently. In addition, we use the model to analyse the impacts of the presence/absence of a FiP.  We examine consumer tariffs, firms' profits, wholesale prices, load shedding, and carbon emissions and we apply our model to a stylised version of the Irish power system for 2030. 

In relation to the research questions outlined in Section \ref{sec:intro}, we can discuss and interpret our findings as follows. \textit{First}, prosumers investing in self-sufficiency can help make markets more efficient. That is, the more consumers that prosume and invest in PV and battery storage, the more we observe the negative effects of market power being mitigated against. While this finding is in line with the literature \citep{devine2023role}, the level of mitigation observed in this work is much greater due to our consideration of investment decisions. However, it should also be noted that, even when all consumers prosume, self-sufficiency is still not enough to fully mitigate against price increases resulting from market power.  

\textit{Second}, we observe that battery storage investments only occur when prices are elevated due to market power. When the market is efficient (perfectly competitive), there is no need for additional storage, particularly not on the demand side, from a systems perspective and for the levels of renewables considered. When market power is assumed, it is only the residential prosumers who invest in storage. They do so in order to protect themselves against the negative effects of market power. Despite this, they would still be better off in a system without market power and thus without the need to invest in storage. 

\textit{Third}, the presence of both market power and self-sufficiency affects the investment decisions of generation firms. The presence of market power increases investment for smaller generation firms but decreases investment for the largest firm we considered. In the absence of prosumers, the assumption of market power increases overall investment in wind and PV by the firms. However, when prosumers are present, market power decreases overall investment by the generation firms. These findings highlight the importance of considering investment decisions on both the demand and supply side in addition to market power. Moreover, the presence of market power alters firms exit decisions, with less decommissioning of conventional peak generation and increased decommissioning of conventional mid-merit generation. This is because market power leads to increased investment in PV and storage on the demand side so that the demand that needs to be covered by the supply-side generators is reduced. 

\textit{Fourth}, self-sufficiency, while constituting and attractive option for consumers/prosumers to reduce their electricity bills, can lead to challenges concerning the cost recovery of the electricity system infrastructure. When continuing to collect these charges (considered as retail premia in this paper) on a volumetric basis and assuming that the gap presented was to be filled with the reduced grid demand, retail premia would need to increase by 10-20\% for industry and 20-40\% for residential consumers. In turn, one would expect that such increased retail premia would further increase the economic attractiveness of self-sufficiency investments, which could easily lead to spiral effect \citep[cf. ][]{mehigan2018review}. Overall, this can result in enormous cost recovery challenges for grid operators and retailers, but also for consumers who cannot invest in self-sufficiency - for whatever reason \citep{bertsch2017drives}. This calls for investigating alternative consumer tariff structures as part of future research.

\textit{Fifth}, we observe that the absence of a FiP, as expected, decreases investment in renewables (particularly wind) by the generation firms. It is interesting to note though that this in turn leads to increased wholesale prices and consequently, increased investment in renewables (particularly PV) by prosumers. Thus, the absence of a FiP leads to  investment in renewables by the prosumers. This finding is counter-intuitive as FiPs are typically introduced to systems in order to incentivize investments in renewables. This also explains why the absence of a FiP only has an attentuated effect on emissions. 

\textit{Sixth}, the presence of market power decreases carbon emissions as market power encourages prosumers to invest in renewables. However, the effects are not as large when all consumers are prosumers. This finding is in contrast to \cite{devine2023role} who, using similar data but do not consider investment decisions, suggest that the presence of market power increases carbon emissions. This finding highlights the importance of considering investment decisions when analysing energy market models. In addition, we find that a FiP reduces carbon emissions, but in the presence of market power and prosumers, the effect of the FiP on emissions becomes relatively small.

As both liberalisation and self-sufficiency are features of modern electricity markets, all of the above findings highlight the importance of including market power in energy market models in addition investment decisions on both the demand side and supply side. Excluding such modelling features may give misleading results. This will be of interest to policymakers, system operators, regulatory, i.e., those that operate on behalf of consumers, as well as generating firms.

Critically reflecting on our approach, we wish to acknowledge some limitations. \textit{First}, modelling all generating firms as Cournot players using a MCP may be unrealistic, particularly given that firm $f=1$ has a generating portfolio larger than all the other firms. As suggested in Section \ref{sec:results}, modelling the largest firm as a Stackelberg leader may be more appropriate. However, the MCP methodology cannot model investment decisions in an oligopoly \citep{devine2023strategic}. Furthermore, MCPs cannot model integer decision variables and thus several aspects of power systems such as online and start-up costs.

Alternative methodologies in the energy market literature include Equilibrium Problem with Equilibrium Constraints (EPECs) and Mathematical Programs with Equilibrium Constraints (MPECs). However, such models can be computationally difficult to solve \citep{pozo2017basic, fanzeres2019robust} and can only be solved for relatively small problems. Hence, EPECs and MPECs would not be appropriate for the current modelling work where we consider a capacity market, 8760 hourly time steps, inter-temporal demand-response constraints, 6 RES-E scenarios, generation firms with multiple technologies and investment decisions, multiple types of consumers who prosume and who may also make investment decisions. The MCP approach outlined in this paper allows us to consider all of these market features together, which is a significant progress beyond the state of the art.

\textit{Second}, each of the players we model is considered to be risk-neutral. In reality, it is more likely that they would make decisions in a risk-averse paradigm \citep{russo2022short}, both in terms of operational and investment decisions. However, modelling risk aversion in equilibrium models is only an emerging area of Operations Research \citep{egging2023stochastic}. Consequently, and given the complexity of the current work, we leave modelling risk-aversion for future work.

\footnotesize
\singlespacing
\bibliography{bibtex_investments}

\end{document}